\documentclass[11pt]{article}
\usepackage{amsmath, graphicx, amsfonts,amssymb, calrsfs}
\usepackage{amsfonts,mathrsfs, color, amsthm}
\usepackage{bm}

\usepackage{enumitem}
\usepackage{hyperref}
\hypersetup{
  colorlinks   = true, 
  urlcolor     = blue, 
  linkcolor    = blue, 
  citecolor   = red 
  }

\def\cK{\mathcal{K}}
\def\sphere{S^{n-1}}

\def\Rn{{\mathbb R^n}}
 \def\R{\mathbb{R}}

\def\dV{\mathcal{V}_{\phi}(K)}
\def\cV{\widetilde{C}_{\phi,\mathcal{V}}}
\def\eV{\mathcal{V}_{\phi}}

\newtheorem{theorem}{Theorem}[section]
\newtheorem{lemma}{Lemma}[section]
\newtheorem{remark}{Remark}[section]
\newtheorem{proposition}{Proposition}[section]
\newtheorem{corollary}{Corollary}[section]

\newtheorem{definition}{Definition}[section]

\def\bt{\begin{theorem}}
\def\et{\end{theorem}}
\def\bl{\begin{lemma}}
\def\el{\end{lemma}}
\def\br{\begin{remark}}
\def\er{\end{remark}}
\def\bc{\begin{corollary}}
\def\ec{\end{corollary}}
\def\bd{\begin{definition}}
\def\ed{\end{definition}}
\def\bp{\begin{proposition}}
\def\ep{\end{proposition}}
\def\ball{B^n_2}

\begin{document}
	\title{On the general dual Orlicz-Minkowski problem }
	\author{Sudan Xing and Deping Ye}
	\date{}
	\maketitle
	\begin{abstract}
	For $K\subseteq \Rn$ a convex body with the origin $o$ in its interior, and $\phi:\mathbb{R}^n\setminus\{o\}\rightarrow(0, \infty)$ a continuous function, define the general dual ($L_{\phi})$ Orlicz  quermassintegral of $K$ by $$\mathcal{V}_\phi(K)=\int_{\mathbb{R}^n \setminus K} \phi(x)\,dx.$$   Under certain conditions on $\phi$, we prove a variational formula for the general dual ($L_{\phi})$ Orlicz quermassintegral, which motivates the definition of $\cV(K, \cdot)$,  the general dual ($L_{\phi})$  Orlicz curvature measure of $K$. 
	
	We pose the following general dual Orlicz-Minkowski problem: {\it Given a nonzero finite Borel measure $\mu$ defined on $\sphere$ and a continuous function $\phi: \Rn\setminus\{o\}\rightarrow (0, \infty)$,  can one find a constant  $\tau>0$ and a convex body $K$ (ideally, containing $o$ in its interior), such that,}  $$\mu=\tau\cV(K,\cdot)?$$ 
Based on the method of  Lagrange multipliers and the established variational formula for the general dual ($L_{\phi})$ Orlicz  quermassintegral,  a solution to the general dual Orlicz-Minkowski problem is provided. In some special cases, the uniqueness of solutions is proved and the solution for $\mu$ being a discrete measure is characterized.

	\end{abstract}

	\section{Introduction} Let $\varphi: (0, \infty)\rightarrow (0, \infty)$ be a continuous function and $\mu$ be a nonzero finite Borel measure defined on the unit sphere $\sphere$, the Orlicz-Minkowski problem asks whether there exists a convex body (a convex and compact subset of $\Rn$ with nonempty interior) $K$ and a constant $\tau>0$, such that, $$\,d\mu=\tau\cdot \varphi(h_K)\,d S_K $$ where $h_K$ denotes the support function of $K$ (see (\ref{support-function--1})) and $S_K$ denotes the surface area measure of $K$ (see (\ref{surface:area:1})).  The Orlicz-Minkowski problem was first investigated by Haberl, Lutwak, Yang and Zhang in their seminal paper \cite{HLYZ2010} for even measure $\mu$. Solutions to the Orlicz-Minkowski problem for $\mu$ being a discrete and/or general (not necessary even) measure were provided by Huang and He \cite{huanghe2012} and Li \cite{liaijun2014}. The planar Orlicz-Minkowski problem in the $L_1$-sense was investigated by Sun and Long \cite{SunLong}. The $p$-capacitary Orlicz-Minkowski problem was posed and studied in \cite{HongYeZhang-1}. The Orlicz-Minkowski problems are central objects in the recent but rapidly developing Orlicz-Brunn-Minkowski theory for convex bodies \cite{GHW2014, Ludwig2010, LYZ2010a,LYZ2010b,XJL}.
	
	The well-studied classical Minkowski problem and its  $L_p$ extension are special cases of the Orlicz-Minkowski problem. When $\varphi(t)=1$, it becomes the classical Minkowski problem back to Minkowski at the turn of the 20th century \cite{min1897,min1903}. Please refer to \cite[Chapter 8]{Sch} for details.  When $\varphi(t)=t^{1-p}$ for $p\in \R$, it becomes the $L_p$ Minkowski problem back to Lutwak \cite{Lu93} in 1993.  Since then, great progress has been made on the $L_p$ Minkowski problem, see e.g., \cite{chen,  chouw06,  HuMaShen, HLW-1, HugLYZ,  JLW-1, Lut-Oli-1, LYZ04, Umanskiy,  zhug2015b, zhug2017}. In particular, the singular cases for $p=0$  and for $p=-n$,  referred to as the logarithmic (or $L_0$) Minkowski problem and the centro-affine Minkowski problem, respectively, are arguably more challenging than the cases for $p\neq 0, -n$. Remarkable contributions on  the logarithmic (or $L_0$) and centro-affine Minkowski problems can be found in, e.g., \cite{BHZ2016, BLYZ2013, chouw06,  JLZ2016, LuWang-1, stancu02, stancu03, stancu08, zhug2014, zhug2015a}. We would like to mention that the $L_p$ Minkowski problem can be described through a fully  nonlinear second-order partial differential equation (i.e., Monge-Amp\`{e}re type equation) and plays fundamental roles in the development of the $L_p$ Brunn-Minkowski theory for convex bodies. 
	
	The $L_p$ surface area measure $h_K^{1-p}S_K$ can be obtained through variational formulas \cite{Lu93, ZHY2017}. As an example, for $p>1$ and $K, L\subseteq \Rn$ convex bodies containing the origin $o$ in their interiors, one has \cite{Lu93} \begin{equation}\label{variational-formula-1} \int_{\sphere}h_L^p(u) h_K^{1-p}(u)\,dS_K(u)=p\cdot \lim_{\varepsilon\rightarrow 0^+} \frac{V(K+_{p} \varepsilon\cdot L)-V(K)}{\varepsilon}\end{equation}  where $V(K)$ stands for the volume of $K$ and $K+_{p} \varepsilon \cdot L$ is a convex body determined by the function $h_{K+_{p}\varepsilon \cdot L}: \sphere \rightarrow (0, \infty)$: for any $\varepsilon>0$, $$h_{K+_{p} \varepsilon\cdot L}^p=h_{K}^p+\varepsilon h_{L}^p.$$ Livshyts  \cite{Livshyts2017} proposed a surface area measure of $K$  with respect to  a measure  $\mu_g$, where $g$, the density of $\mu_g$ with respect to the Lebesgue measure, is  continuous on its support. A variational formula for $\mu_g$ similar to (\ref{variational-formula-1}) for $p=1$ was also provided in \cite{Livshyts2017}, which gives a variational interpretation of the surface area measure of $K$  with respect to $\mu_g$. The related Minkowski problem was posed and a solution to this problem was given under certain conditions on $\mu_g$ (such as, $\mu_g$ being a measure with positive degree of concavity and positive degree of homogeneity). An $L_p$ extension of the theory by Livshyts was obtained by Wu \cite{Wu}, where the $L_p$ surface area measure with respect to $\mu_g$ was proposed and related $L_p$ Minkowski problem was solved under certain conditions on $\mu_g$.  Indeed, our paper was partially motivated by  \cite{Livshyts2017, Wu}. 
	
	This paper is also motivated by the recent work of Zhu, Xing and Ye on the dual Orlicz-Minkowski problem \cite{ZSY2017}, which belongs to the recently initiated dual Orlicz-Brunn-Minkowski theory \cite{ghwy15, Ye2016a, Zhub2014} and can be viewed as the ``dual" of the Orlicz-Minkowski problem. For a convex body $K\subseteq \Rn$ with the origin $o$ in its interior and a continuous function $\varphi: (0, \infty)\rightarrow (0, \infty)$, the authors in \cite{ZSY2017} defined   the dual Orlicz curvature measure $\widetilde{C}_\varphi(K,\cdot)$ and investigated the following  dual Orlicz-Minkowski problem: {\it under what conditions on $\varphi$ and a given nonzero finite Borel measure $\mu$ on $\sphere$,  there exist a constant $\tau>0$ and a convex body $K$ (ideally with the origin in its interior) such that $\mu=\tau \widetilde{C}_\varphi(K,\cdot)$?} A solution to the dual Orlicz-Minkowski problem was given under the assumptions: a) the measure $\mu$ is not concentrated on any closed hemishpere, i.e., $\mu$ satisfying (\ref{condition for Minkowski problem}); b) the function $\varphi$ and its companion function \begin{equation*}  \phi(t)=\int_t^\infty \frac{\varphi(s)}{s}\,ds\end{equation*} satisfy conditions A1)-A3) as described in Section \ref{Section:2}. We would like to mention that the assumption on $\mu$, i.e., (\ref{condition for Minkowski problem}), is necessary for the solutions of various Minkowski problems. A special case with $\varphi(t)=t^q$ for $q<0$ was solved in the remarkable paper \cite{zhao} by Zhao, as $\varphi(t)=t^q$ and its companion $\phi(t)=-t^q/q$ satisfy conditions A1)-A3). The dual Orlicz-Minkowski problem stemmed from the groundbreaking work \cite{HLYZ} in 2016 by Huang, Lutwak,  Yang and Zhang, where they provided a very detailed study of the geometric measures (such as the $q$-th dual curvature measures) in the dual Brunn-Minkowski theory and initiated the very promising dual Minkowski problem for the $q$-th dual curvature measures. In particular, they provided a solution to the dual Minkowski problem for the $q$-th dual curvature measures with $q\in (0, n]$  and even measure $\mu$ (plus some additional conditions). Note that the logarithmic Minkowski problem is the case for $q=n$. Since their groundbreaking work \cite{HLYZ}, there is a growing body of work in this direction, see e.g., \cite{BHP, BLYZZ2017, Henk,  JiangWu, LiSW-1,  zhao, zhao-jdg, ZSY2017}. 
	
	The starting point of this paper is the general dual ($L_{\phi})$ Orlicz quermassintegral. For $K\subseteq \Rn$ a convex body with the origin $o$ in its interior, and $\phi:\mathbb{R}^n\setminus\{o\}\rightarrow(0, \infty)$ a continuous function, define the general dual ($L_{\phi})$ Orlicz quermassintegral of $K$  by $$\mathcal{V}_\phi(K)=\int_{\mathbb{R}^n \setminus K} \phi(x)\,dx.$$ In order to have  $\mathcal{V}_\phi(K)$ well-defined for each convex body $K$ with the origin $o$ in its interior and to solve the general dual Orlicz-Minkowski problem, some conditions on $\phi$ are required and these conditions are described in Section \ref{Section:2} (i.e.,  conditions C1) and C2) following Definition \ref{the general dual orlicz quermassintegral}). Some special functions satisfying conditions C1) and C2) are discussed. The convergence of $\eV$ is summarized in Lemma \ref{continuity-general-dual-qu-1}, and it will be used to establish the existence of the solutions to the general dual Orlicz-Minkowski problem. 
	
	To formulate the general dual Orlicz-Minkowski problem, the general dual ($L_{\phi})$  Orlicz curvature measure is required. For $K\subseteq \Rn$ being a convex body with the origin $o$ in its interior and a subset $\eta\subseteq \sphere$, denote by $\rho_K$ the radial function of $K$ and $\pmb{\alpha}^*_K(\eta)\subseteq \sphere$  the reverse radial Gauss image of $\eta$, respectively. Define $\widetilde{C}_{\phi,\mathcal{V}}(K,\cdot)$,  the general dual ($L_{\phi})$ Orlicz curvature measure of $K$ with $\phi :\mathbb{R}^n\setminus\{o\}\rightarrow(0, \infty)$  a continuous function satisfying condition C1) (following Definition \ref{the general dual orlicz quermassintegral} in Section \ref{Section:2}), by 
	$$ \widetilde{C}_{\phi,\mathcal{V}}(K,\eta)=\int_{\pmb{\alpha}^*_K(\eta)}\phi(\rho_K(u)u)[\rho_K(u)]^n \,du$$  for any Borel set $\eta\subseteq \sphere$, where $\,du$ is the spherical measure of $\sphere$.  The properties for the general dual ($L_{\phi})$  Orlicz curvature measure are provided in Section \ref{Section:3-1}. In particular, convenient formulas to calculate integrals with respect to $\widetilde{C}_{\phi,\mathcal{V}}(K,\cdot)$ are given in Lemma \ref{two formula for dual orlicz curvature measure}, and the weak convergence of the general dual ($L_{\phi})$  Orlicz curvature measure is summarized in Proposition \ref{weak continuity}. These properties are crucial in solving the general dual Orlicz-Minkowski problem:  {\it Given a nonzero finite Borel measure $\mu$ defined on $\sphere$ and a continuous function $\phi: \Rn\setminus\{o\}\rightarrow (0, \infty)$,  can one find a constant  $\tau>0$ and a convex body $K$ (ideally, containing $o$ in its interior), such that,}  $$\mu=\tau\cV(K,\cdot)?$$

 A basic method to solve various Minkowski problems is the method of  Lagrange multipliers, and hence a variational formula related to the general dual ($L_{\phi})$  Orlicz curvature measure is essential.  Such a variational formula is proved in Theorem \ref{variational-for-decreaing-1}. In fact, this variational formula motivates our definition of $\cV(K, \cdot)$ and provides a variational interpretation for the general dual ($L_{\phi})$  Orlicz curvature measure.  With the help of the method of  Lagrange multipliers and the established variational formula in Theorem  \ref{variational-for-decreaing-1},  a solution to the general dual Orlicz-Minkowski problem is provided in Theorem \ref{solution-general-dual-Orlicz-main theorem}. Conditions C1) and C2)  are listed in Section \ref{Section:2} following Definition \ref{the general dual orlicz quermassintegral}. 
 
 \vskip 2mm \noindent {\bf Theorem \ref{solution-general-dual-Orlicz-main theorem}.} {\em  Let $\mu$ be a nonzero finite Borel measure on $S^{n-1}$ satisfying (\ref{condition for Minkowski problem}) and let $\phi$ be a function satisfy conditions C1) and C2).
	Then there exists a convex body $K$ containing the origin $o$ in its interior, such that, $$\frac{\mu}{|\mu|}=\frac{\widetilde{C}_{\phi,\mathcal{V}}(K,\cdot)}{\widetilde{C}_{\phi,\mathcal{V}}(K,S^{n-1})}, $$ where $|\mu|=\mu(\sphere)$ is the total $\mu$-mass of $\sphere$.}
	
	 It seems, in general, not possible to obtain the uniqueness of the solutions to the  general dual Orlicz-Minkowski problem. However, in some special cases, say $\phi$ having certain homogeneity, we are able to prove the uniqueness. Moreover,  the solution to the general dual Orlicz-Minkowski problem is proved to be a polytope, when $\mu$ is a discrete measure. The details will be provided in Section \ref{Section:6}.

	\section{The general dual Orlicz quermassintegral}\label{Section:2}
Our setting in this article is the $n$-dimensional Euclidean space $\mathbb{R}^n$ equipped with the standard Euclidean norm $|\cdot|$ induced by the inner product $\langle\cdot, \cdot \rangle$.  The standard notations $o$, $\ball$ and $\sphere$ denote the origin, the unit Euclidean ball and the unit sphere in $\mathbb{R}^n$, respectively.  For a set $K\subseteq \R^n$,  the boundary of $K$ and the interior of $K$ are denoted by $\partial K$ and $\mathrm{int} K$, respectively. The unit vector $\bar{x}=x/|x|\in \sphere$ refers to the direction vector of $x\in \Rn\setminus\{o\}$. By $\mathrm{conv}(A)$, we mean the convex hull of $A\subset \Rn$; namely, $\mathrm{conv}(A)$ is the smallest convex set containing $A$. 

We consider the measure $\mathcal{V}_\phi$, whose density function with respect to the Lebesgue measure $\,dx$ is a continuous function $\phi:\mathbb{R}^n\setminus\{o\}\rightarrow(0, \infty)$.   \bd\label{the general dual orlicz quermassintegral} For a measurable subset  $E\subseteq \Rn$ with $o\in \mathrm{int} E$, define $\mathcal{V}_\phi(E)$ by  $$\mathcal{V}_\phi(E)=\int_{\mathbb{R}^n \setminus E} \phi(x)\,dx.$$\ed Clearly, $\eV(\cdot)$ is monotone decreasing, that is, if $E\subseteq F$ with $o\in \mathrm{int} E$, then $\Rn\setminus E\supseteq \Rn\setminus F$ and hence $$\eV(E)=\int_{\mathbb{R}^n \setminus E} \phi(x)\,dx\geq \int_{\mathbb{R}^n \setminus F} \phi(x)\,dx=\eV(F)$$ due to the positivity of $\phi$. 
	
When $E$ is a star-shaped set in $\Rn$, $\mathcal{V}_\phi(E)$ can be reformulated through the radial function of $E$ and the spherical measure $\,du$ on $\sphere$. Hereafter, $E\subseteq \Rn$ is said to be a star-shaped set with  respect to $o$, if $o\in E$  and the line segment $[o,x]\subseteq E$ for all $x\in E$.   For a star-shaped set $E$ with  respect to $o$, one can define its radial function  $\rho_E: \sphere \rightarrow [0, \infty]$  by   $$
	\rho_E(u)=\sup\big\{\lambda>0:\ \lambda u\in E\big\} \ \ \ ~\text{for\ each}\ u\in\sphere.$$ Denote by $\mathcal{S}$ the set of all star-shaped sets in $\R^n$ with respect to $o$ whose radial functions are measurable.  In fact, we will be working on the set $\mathcal{K}_{o}^n\subseteq \mathcal{S}$, the collection of all convex bodies in $\R^n$ containing the origin $o$ in their interiors. That is, if $K\in \mathcal{K}_{o}^n$, then $K\subseteq \Rn$ is a convex compact set with the origin $o$ in its interior.   The support function of $K\in\mathcal{K}_{o}^n$, $h_{K}(u): \sphere\rightarrow\mathbb{R}$,  is defined by
	\begin{equation}\label{support-function--1}
	h_{K}(u)=\sup_{x\in K}\langle x, u\rangle\ \ \ \ ~\text{for\  each}\ u\in\sphere.
	\end{equation} 
When $E\in \mathcal{S}$, $\eV(E)$ can be calculated by \begin{eqnarray} \mathcal{V}_\phi(E) = \int_{\mathbb{R}^n \setminus E}\phi(x)\,dx = \int_{\sphere} \bigg(\int_{\rho_E(u)}^\infty \phi(ru)r^{n-1}\,dr\bigg)\,du.\label{general-quermass-11}
	\end{eqnarray}  For convenience, let $$\Phi(t,u)=\int_{t}^{\infty}\phi(ru)r^{n-1}\,dr$$ and hence formula (\ref{general-quermass-11}) can be rewritten as \begin{eqnarray}\mathcal{V}_\phi(E)=\int_{\sphere}\Phi(\rho_{E}(u),u)\,du. \label{general-quermass-22}
	\end{eqnarray}

  For each $K\in\mathcal{K}_{o}^n$, one can define $K^*$,  the polar body of $K$,  by
	\[
	K^*=\{x\in\mathbb{R}^n: \langle x, y\rangle\leq1 \ \text{for all}\  y\in K \}.\]  Clearly $K^*\in \mathcal{K}_{o}^n$. Moreover, the bipolar theorem asserts that $(K^*)^*=K$ (see e.g., \cite{Sch}) and then $$\rho_K (u) h_{K^*}(u)=h_K (u) \rho_{K^*}(u)=1.$$ Hence, for each $K\in\mathcal{K}_{o}^n$, one gets, by formula  (\ref{general-quermass-22})
 \begin{eqnarray}\mathcal{V}_\phi(K)=\int_{\sphere}\Phi(\rho_{K}(u),u)\,du=\int_{\sphere}\Phi(h_{K^*}(u)^{-1},u)\,du. \label{general-quermass-33}
	\end{eqnarray} In later context,  for each $K\in \mathcal{K}_{o}^n$,  $\mathcal{V}_\phi(K)$ will be called the general dual ($L_{\phi})$ Orlicz quermassintegral of $K$.  
	
Now we list the basic conditions for function $\phi$:  \begin{itemize} \item[C1)]  $\phi:\mathbb{R}^n\setminus\{o\}\rightarrow(0, \infty)$ is a continuous function, such that,  for any fixed $t>0$, the function 
 $$ \Phi(t, u)=\int_{t}^{\infty}\phi(ru)r^{n-1}\,dr$$ is positive and continuous on $S^{n-1}$;  \item [C2)] for any fixed $u_{0} \in \sphere$ and any fixed positive constant $b_{0}\in (0, 1)$, one has   $$\lim_{a\rightarrow 0^+}\eV(\mathbb{R}^n\setminus \mathcal{C}(u_{0},a,b_{0}))=\infty, $$ where $\mathcal{C}(u_{0},a,b_{0})$ is defined by $$\mathcal{C}(u_{0},a,b_{0})=\Big\{x\in\mathbb{R}^n:\ \langle \bar{x}, u_{0}\rangle\geq b_{0}\ \mathrm{and} \ |x|\geq a\Big\}.$$ 
	\end{itemize}
	 	 	
		In fact, condition C1)  guarantees that $\eV(K)<\infty$ for each $K\in \mathcal{K}_o^n$. To see this, as $o\in \mathrm{int} K$, there exists a constant $r_0>0$ such that $ r_0\ball \subseteq K$. By formula (\ref{general-quermass-33})  and the fact that $\eV(\cdot)$ is monotone decreasing, one has, 
		\begin{eqnarray*}\mathcal{V}_\phi(K)\leq \eV(r_0\ball)=\int_{\sphere}\Phi(r_0,u) \,du<\infty. 
	\end{eqnarray*} Condition C2) is for the solution of the general dual Orlicz-Minkowski problem.

	A typical function satisfying conditions C1) and C2) is a continuous function  $\phi: \Rn\setminus\{o\}\rightarrow (0, \infty)$ such that 
	\begin{equation}\label{inf-sup-condition} \sup_{|x|>r_1}\phi(x)|x|^{n-\alpha_1-1}\leq C_1 \ \ \ \mathrm{and} \ \ \ \inf_{|x|<r_1}\phi(x)|x|^{n-\alpha_2-1}\geq C_2 \end{equation}  hold for some constants $0<r_1<\infty$, $C_1<\infty$, $C_2>0$ and $-\infty<\alpha_1, \alpha_2<-1$. In particular, if   
	$$\lim_{|x|\rightarrow\infty}\phi(x)|x|^{n-\alpha_1-1}=C_1 \ \ \mathrm{and} \ \ \lim_{|x|\rightarrow 0}\phi(x)|x|^{n-\alpha_2-1}=C_2$$ for some constants $0< C_1, C_2<\infty$ and $-\infty<\alpha_1, \alpha_2<-1$, then such $\phi$ satisfies (\ref{inf-sup-condition}) (for different constants). Now let us check that a continuous function $\phi$ satisfying (\ref{inf-sup-condition}) must also satisfy conditions C1) and C2). To this end, let $t>0$ and  $u\in S^{n-1}$ be fixed. It is obvious to have $\Phi(t, u)>0$. Moreover  
		 \begin{eqnarray*}
			\Phi(t, u) &= & \int_{t}^{r_{1}}\phi(ru)r^{n-1}\,dr +\int_{r_{1}}^{\infty}\phi(ru)r^{n-1}\,dr\\ &\leq & \Big|\int_{t}^{r_{1}}\phi(ru)r^{n-1}\,dr\Big|+\int_{r_{1}}^{\infty}\phi(ru)r^{n-1}\,dr\\
			&\leq& \Big|\int_{t}^{r_{1}}\phi(ru)r^{n-1}\,dr\Big| +C_{1} \int_{r_{1}}^{\infty}r^{\alpha_{1}}\,dr \\
			&=&\Big|\int_{t}^{r_{1}}\phi(ru)r^{n-1}\,dr\Big|-\frac{C_{1}}{\alpha_{1}+1}\cdot r_{1}^{\alpha_{1}+1}.						
		\end{eqnarray*} Thus $\Phi(t, u) <\infty$ due to the continuity of $\phi$, and $\Phi(t, u)$ is well defined. Now we claim that $\Phi(t, \cdot)$ is continuous on $\sphere$.  For fixed $t$ and for an arbitrary sequence $u_i\rightarrow u$ with $u_i, u\in \sphere$, one has, for all $r\geq t$,  $\phi(ru_{i})r^{n-1}\rightarrow \phi(ru)r^{n-1}$ and $$\phi(ru_{i})r^{n-1}\leq C_1 r^{\alpha_1} +M$$  for all $i\geq 1$, where, due to the continuity of $\phi$,  $$M=\max\Big\{\phi(x)|x|^{n-1}: \ \ |x| \ \mathrm{is\ between\ } t \ \mathrm{and} \  r_1 \Big\}<\infty.$$ It follows from the dominated convergence theorem that  \[
\lim_{i\rightarrow \infty} \Phi(t,u_{i})=\lim_{i\rightarrow \infty}  \int_{t}^{\infty}\phi(ru_{i})r^{n-1}\,dr=  \int_{t}^{\infty} \lim_{i\rightarrow \infty} \phi(ru_{i})r^{n-1}\,dr = \int_{t}^{\infty}\phi(ru)r^{n-1}\,dr=\Phi(t,u).
\] Hence   $\Phi(t,u)$ is continuous on $\sphere$ and   C1) is verified. Now let us verify C2)  as follows: for any $b_0\in (0, 1)$,	\begin{eqnarray*}
	\lim_{a\rightarrow 0^+}\eV(\mathbb{R}^n\setminus \mathcal{C}(u_{0},a,b_{0})) &=&\lim_{a\rightarrow 0^+}	
	\int_{\{u\in \sphere: \langle u, u_0\rangle\geq b_0\}}\int_{a}^{\infty}\phi(ru)r^{n-1}\,dr\,du\\&\geq&\limsup_{a\rightarrow 0^+}	
	\int_{\{u\in \sphere: \langle u, u_0\rangle\geq b_0\}}\int_{a}^{r_1}\phi(ru)r^{n-1}\,dr\,du\\	&\geq&	C_{2}\cdot \limsup_{a\rightarrow 0^+}
	\int_{\{u\in \sphere: \langle u, u_0\rangle\geq b_0\}}  \int_{a}^{r_{1}}r^{\alpha_{2}}\,dr\,du\\
	&=&C_{2}\cdot \bigg(\int_{\{u\in \sphere: \langle u, u_0\rangle\geq b_0\}} \,du\bigg)  \cdot \limsup_{a\rightarrow 0^+} \frac{r_1^{1+\alpha_2}-a^{1+\alpha_2}}{1+\alpha_2}\\ &=& \infty,
\end{eqnarray*} where we have used (\ref{general-quermass-11}), (\ref{inf-sup-condition}) and  $\alpha_2<-1$.

Now let us provide several special cases of functions satisfying conditions C1) and C2). 
 
\vskip 2mm \noindent {\it Case 1:} $\phi(x)=\psi(|x|)$ for all $x\in \Rn\setminus \{o\}$ with $\psi: (0, \infty)\rightarrow (0, \infty)$ a continuous function. In this case, \begin{equation}\label{definition:psi:hat} \Phi(t, u)=  \int_{t}^{\infty}\phi(ru)r^{n-1}\,dr=  \int_{t}^{\infty}\psi(r)r^{n-1}\,dr:=\frac{1}{n}  \cdot \hat{\phi}(t).\end{equation} Equivalently, \begin{equation}\label{def:hat:psi-phi} \psi(t)=- \hat{\phi}'(t)t^{1-n}/n. \end{equation}    By formula (\ref{general-quermass-33}), one has,  for $K\in \cK_o^n$,
	\[
	\eV(K)=\int_{\sphere}\Phi(\rho_{K}(u), u)\,du=\frac{1}{n}  \int_{\sphere}\hat{
		\phi}(\rho_{K}(u))\,du=\widetilde V_{\hat{
			\phi}}(K),
	\] where $\widetilde V_{\varphi}(\cdot)$  is the dual ($L_{\varphi})$ Orlicz quermassintegral proposed in \cite{ZSY2017}, namely  
	\begin{equation*} 
		\widetilde{V}_\varphi(K)= \frac{1}{n}\int_{\sphere}\varphi(\rho_K(u))du.
	\end{equation*}  In  \cite{ZSY2017}, the dual Orlicz-Minkowski problem is solved under the following conditions: 
	\begin{itemize} \item[A1)]  $\hat{
			\phi}: (0, \infty)\rightarrow (0, \infty)$ is a strictly decreasing continuous function with $$\lim_{t\rightarrow  0^+}\hat{
			\phi}(t)=\infty \ \ \ \mathrm{and} \ \ \ \lim_{t\rightarrow  \infty}\hat{\phi}(t)=0;$$
		\item [A2)]   $\hat{
			\phi}'$, the derivative of $\hat{
			\phi}$, exists and is strictly negative on $(0, \infty)$; \item [A3)]  $\hat{\varphi}(t)=-\hat{
			\phi}'(t)t: (0, \infty)\rightarrow (0, \infty)$ is continuous; hence    \begin{equation*}  \hat{
				\phi}(t)=\int_t^\infty \frac{\hat{\varphi}(s)}{s}\,ds. \end{equation*} \end{itemize} In  Case 1,  it is obvious that $\hat{\varphi}(t)=n\psi(t)t^n$.  Now let us check that if $\hat{\phi}$ and its companion function $\hat{\varphi}$ satisfy conditions   A1)-A3), then $\phi(x)=\psi(|x|)$ with $\psi$ given by (\ref{def:hat:psi-phi})  satisfies conditions C1) and C2). In fact, condition C1) can be easily checked by (\ref{definition:psi:hat}) and A1). Let us verify condition C2) as follows: for any $b_0\in (0, 1)$, \begin{eqnarray*}
	\lim_{a\rightarrow 0^+}\eV\big(\mathbb{R}^n\setminus \mathcal{C}(u_{0},a,b_{0})\big) &=&\lim_{a\rightarrow 0^+}	
	\int_{\{u\in \sphere:\langle u, u_{0}\rangle \geq b_{0}\}}\int_{a}^{\infty}\phi(ru)r^{n-1}\,dr\,du\\
	 &=& \frac{1}{n}\cdot  \lim_{a\rightarrow 0^+}	\hat{\phi}(a) \cdot \bigg(
	\int_{\{u\in \sphere:\langle u, u_{0}\rangle \geq b_{0}\}}\,du\bigg)\\
	&=&\infty,	
\end{eqnarray*}	  where we have used (\ref{general-quermass-11}), (\ref{definition:psi:hat}), and  condition A1).

	  \vskip 2mm \noindent {\it Case 2:}   $\phi(x)=\psi(|x|)\phi_{2}(\bar{x})$ where $\bar{x}=x/|x|$,  $\psi: (0, \infty)\rightarrow (0, \infty)$ is a continuous function  on $(0, \infty)$, and $ \phi_{2}: \sphere \rightarrow (0, \infty)$ is a continuous function on $\sphere$.   In this case,  the general dual ($L_{\phi})$ Orlicz quermassintegral of $K\in \mathcal{K}_o^n$  has the following form: \begin{eqnarray} \mathcal{V}_\phi(K)&=&\int_{\sphere} \int_{\rho_K(u)}^{\infty}\phi(ru)r^{n-1}\,dr \,du \nonumber \\ &=& 
	  \int_{\sphere}  \bigg( \int_{\rho_K(u)}^{\infty}\psi(r)r^{n-1}\,dr \bigg)  
	\phi_2(u)  \,du \nonumber \\ &=&\frac{1}{n}\int_{S^{n-1}}\hat{\phi}(\rho_K(u))\phi_{2}(u)\,du,\label{homogenes-1-1}
	\end{eqnarray} 
	where $\hat{\phi}$ is given by (\ref{definition:psi:hat}). Again, if $\hat{\phi}$ and its companion function $\hat{\varphi}$ satisfy conditions   A1)-A3), then  $\phi(x)=\psi(|x|)\phi_{2}(\bar{x})$ with $\psi$ give by (\ref{def:hat:psi-phi}) satisfies conditions C1) and C2); this follows from an argument similar to the one as in Case 1. A typical example in this case is $$\phi(x)=\|x\|^{q-n}=|x|^{q-n}\cdot \|\bar{x}\|^{q-n}$$ where $q<0$ is a constant and $\|\cdot\|: \Rn\rightarrow [0, \infty)$ is any norm on $\Rn$. (Note that $\phi_2(\bar{x})= \|\bar{x}\|^{q-n}$ is always positive, due to the equivalence between the two norms $\|\cdot\|$ and $|\cdot|$.)  Indeed, when $\phi(x)=\|x\|^{q-n}=|x|^{q-n}\cdot \|\bar{x}\|^{q-n}$, then $\psi(|x|)=|x|^{q-n}$. Hence $$\hat{\phi}(t)=n\int_{t}^{\infty}\psi(r)r^{n-1}dr=n\int_{t}^{\infty}r^{q-1}\,dr=-\frac{n}{q} \cdot t^q$$ and  $\hat{
	\varphi}=nt^{q},$ which satisfy conditions A1)-A3).

  A sequence of convex bodies $\{K_i\}_{i=1}^{\infty}\subseteq \mathcal{K}_{o}^n$ converging to a convex body  $K\in \mathcal{K}_o^n$ in the sense of Hausdorff metric means that
	\begin{equation} \label{H-metric-1}
	\|h_{K_i}-h_K\|_{\infty}=\sup_{u\in\sphere} |h_{K_i}(u)-h_K(u)|\rightarrow0\ \ \  \mathrm{as} \ \ \ i\rightarrow\infty.
	\end{equation}  Indeed, this is equivalent to \begin{eqnarray}	\|\rho_{K_i}-\rho_K\|_{\infty}=\sup_{u\in\sphere} |\rho_{K_i}(u)-\rho_K(u) |\rightarrow0\ \ \  \mathrm{as} \ \ \ i\rightarrow\infty. \label{conv-Hausdorff-radial}
	\end{eqnarray}  We will need the following convergence result regarding $\eV(\cdot)$. 

	\bl\label{continuity-general-dual-qu-1} Assume that  $\phi$ is  a function  satisfying condition C1). If the sequence $\{K_i\}_{i=1}^{\infty} \subseteq \cK_o^n$ converges to $K\in\cK_o^n$ in the sense of Hausdorff metric, then
	$$\lim_{i\rightarrow \infty}\mathcal{V}_\phi(K_i)=\dV.$$
	\el
	\begin{proof}   Let $\phi:\R^n \setminus\{o\} \rightarrow (0, \infty)$ be  a continuous function satisfying C1).  It can be checked that,  for any fixed $u\in\sphere$  and for any fixed constant  $t_0>0$,   	\begin{equation}\label{continuity for general dual orlicz quermassintegral 1}
		\lim_{t\rightarrow t_0}\Phi(t,u)=\Phi(t_0,u).
	\end{equation}  In fact, 	for any fixed $u\in\sphere$, $\Phi(t,u)$ is a decreasing function on $t\in (0, \infty)$. Let $t\rightarrow t_0$, and without loss of generality assume that $t>t_0/2$. By condition C1) and the fact that $\Phi(t,u)$ is decreasing on $t$, one has, $$\Phi(t,u)=\int_{t}^{\infty}\phi(ru)r^{n-1}\,dr\leq \int_{t_0/2}^{\infty}\phi(ru)r^{n-1}dr=\Phi(t_0/2, u)<\infty.$$ It follows from the dominated convergence theorem that  \begin{eqnarray*} \lim_{t\rightarrow t_0}\Phi(t,u)= \lim_{t\rightarrow t_0}  \int_{t}^{\infty}\phi(ru)r^{n-1}\,dr =  \int_{t_0}^{\infty}\phi(ru)r^{n-1}\,dr =\Phi(t_0,u).
	\end{eqnarray*}  

Let $\{K_i\}_{i=1}^{\infty}\subseteq \mathcal{K}_o^n$ be a sequence of convex bodies converging to $K\in \mathcal{K}_o^n$ in the Hausdorff  metric. Based on (\ref{conv-Hausdorff-radial}), $\rho_{K_i}$ converges to $\rho_K$ uniformly on $\sphere$. Moreover, as $K\in \mathcal{K}_o^n$, one can find a constant $R_{1}>0$, such that, for all $u\in \sphere$ and for all $i=1, 2, \cdots$, $$R_{1}\ \le \ \rho_{K_i}(u)  \ \ \ \mathrm{and} \ \ \ R_1\leq \rho_{K}(u).$$
Together with the fact that $\Phi(t,u)$ is a decreasing function on $t\in (0, \infty)$, one has    
		\[\Phi(\rho_{K_{i}}(u),u) \leq\Phi(R_{1},u) \ \ \mathrm{and} \ \ \Phi(\rho_{K}(u),u)\leq\Phi(R_{1}, u) \ \ \ \text{for\ all}\ u\in\sphere. 
		\] By condition C1), $\Phi(R_{1}, u)$ is positive and continuous on $S^{n-1}$. Hence, $$\int_{\sphere} \Phi(R_{1}, u)\,du<\infty.$$ It follows from (\ref{general-quermass-33}), 
		(\ref{continuity for general dual orlicz quermassintegral 1}) and the dominated convergence theorem that
		\begin{eqnarray*}
			\lim_{i\rightarrow\infty}\mathcal{V}_\phi(K_i)&=&\lim_{i\rightarrow\infty} \int_{\sphere}\Phi(\rho_{K_{i}}(u),u)\,du\\ &=&\int_{\sphere}\lim_{i\rightarrow\infty}\Phi(\rho_{K_{i}}(u),u)\,du\\&=&\int_{\sphere}\Phi(\rho_{K}(u),u)\,du\\ &=&\eV(K). \end{eqnarray*} This concludes the proof of Lemma \ref{continuity-general-dual-qu-1}. \end{proof}
  	
	\section{The general dual Orlicz curvature measure}\label{Section:3-1}
  For $K\in \mathcal{K}_o^n$,  the supporting hyperplane of $K$ at the direction $u\in \sphere$, denoted by $H(K,u)$, is given by $$
H(K,u)=\big\{x\in\R^n: \langle x, u\rangle=h_K(u)\big\}.$$ Denoted by $\pmb{\alpha}^*_K(\eta)$ the reverse radial Gauss image of $\eta\subseteq \sphere$, that is, $$
\pmb{\alpha}^*_K(\eta)=\big\{\bar{x}=x/|x|: x\in \partial K \cap H(K,u)\ \ \text{for\ some}\ u\in \eta\big\}.
$$  For convenience, let $$\Psi_{K}(u)=\phi(\rho_K(u)u)[\rho_K(u)]^n \ \ \ \mathrm{for}\ u\in\sphere.$$ In fact,  for any $x\in \partial K$, one has $\Psi_{K}(\bar{x})=\phi(x)|x|^n.$   

We are ready to give the definition  of the general dual Orlicz curvature measure.
	\bd \label{general dual Orlicz curvature measure}
	For any $K\in \mathcal{K}_o^n$ and for any function $\phi$ satisfying condition C1),  the general dual ($L_{\phi})$ Orlicz curvature measure of $K$, denoted by $\widetilde{C}_{\phi,\mathcal{V}}(K,\cdot)$, is given by 
	$$
	\widetilde{C}_{\phi,\mathcal{V}}(K,\eta)=\int_{\pmb{\alpha}^*_K(\eta)}\Psi_{K}(u)\,du$$  for any Borel set $\eta\subseteq \sphere$. 
	\ed
	Indeed, for each $K\in \mathcal{K}_o^n$,  $\widetilde{C}_{\phi,\mathcal{V}}(K,\cdot)$ does define a  Borel measure  on $\sphere$. 
	 To this end, we only need to show that $\cV(K,\cdot)$ satisfies the countable additivity, as $\cV(K,\emptyset)=0$ holds trivially. That is, we need to prove \begin{eqnarray*}\cV(K,\cup_{i=1}^{\infty}\eta_{i})=\sum_{i=1}^\infty\cV(K,\eta_i)
	\end{eqnarray*} for any sequence of pairwise disjoint Borel sets  $\eta_1, \eta_2, \cdots \subseteq S^{n-1}$. Recall that $
	\pmb{\alpha}^*_K\big(\!\cup_{i=1}^\infty\eta_i\big)
	=\cup_{i=1}^\infty\pmb{\alpha}^*_K(\eta_i)
	$ by \cite[Lemma 2.3]{HLYZ} and  $$\pmb{\alpha}_{K}^*(\eta_i)=\overline{\pmb{x}_K(\eta_i)}=\{\bar{x}: \ x\in \pmb{x}_K(\eta_i)\}\subseteq\sphere$$   is spherical measurable for each $i\geq 1$ by \cite[Lemma 2.1]{HLYZ}, where $\pmb{x}_K(\eta_i)$ is the reverse spherical image of  $\eta_i\subseteq\sphere$ given by   
	$$\pmb{x}_K(\eta_i)=\big\{x\in\partial K: x\in H(K,u)\ \ \text{for some}\ \ u\in \eta_i\big\}\subseteq \partial K.$$ Therefore,
	\begin{eqnarray}
		\cV(K,\cup_{i=1}^\infty\eta_i)
		=\int_{\pmb{\alpha}^*_K(\cup_{i=1}^\infty\eta_i)}\Psi_{K}(u)\,du
		=\int_{\cup_{i=1}^\infty\pmb{\alpha}^*_K(\eta_i)}\Psi_{K}(u)\,du. \label{additivity--111}
	\end{eqnarray} The  countable additivity will follow immediately if $\cup_{i=1}^\infty\pmb{\alpha}^*_K(\eta_i)$ is pairwise disjoint.  However,   by \cite[Lemma 2.4]{HLYZ}, one gets that $\{\pmb{\alpha}_K^*(\eta_j)\setminus \omega_K\}_{j=1}^\infty$ is pairwise disjoint, where   	\[
	\omega_K= \big\{v \in \sphere: \  \pmb{\alpha}_K(v) \ \ \mathrm{has\ more\ than\ one\ element}\big\}
	\] with $\pmb{\alpha}_K(\omega)$, the radial Gauss image  of $\omega\subseteq S^{n-1}$, given by 
	$$
	\pmb{\alpha}_K(\omega)=\pmb{\nu}_K(\{\rho_K(u)u\in\partial K: u\in \omega\})
	\subseteq \sphere$$ and with $\pmb{\nu}_K(\sigma)$, the spherical image of $\sigma \subseteq \partial K$,  given by
	$$
	\pmb{\nu}_K(\sigma)=\big\{u\in S^{n-1}: x\in H(K,u)\ \ \text{for some}\ \ x\in \sigma\big\}\subseteq \sphere. 
	$$ Fortunately, the spherical
	measure of $\omega_K$ turns out to be zero \cite[p.339-340]{HLYZ} and hence, by (\ref{additivity--111}),  \begin{eqnarray*}
		\cV(K,\cup_{i=1}^{\infty}\eta_{i})  &=& \int_{\cup_{i=1}^\infty(\pmb{\alpha}^*_K(\eta_i)\setminus \omega_K)}\Psi_{K}(u)\,du\\
		&=& \sum_{i=1}^\infty\int_{\pmb{\alpha}^*_K(\eta_i)\setminus \omega_K}\Psi_{K}(u)du\\
		&=& \sum_{i=1}^\infty\int_{\pmb{\alpha}^*_K(\eta_i)}\Psi_{K}(u)du\\
		&=&\sum_{i=1}^\infty\cV(K,\eta_i).
	\end{eqnarray*}
	This concludes that $\cV$ is a Borel measure.

	Note that $\pmb{\nu}_K(x)$ for $x\in \partial K$, $\pmb{x}_K(u)$ and $\pmb{\alpha}_K(u)$ for $u\in \sphere$ may contain more than one element. When they are singleton sets, $\nu_K(x)$ for $x\in \partial K$, $x_K(u)$ and $\alpha_K(u)$ for $u\in \sphere$ are used, respectively,  instead of $\pmb{\nu}_K(x)$, $\pmb{x}_K(u)$ and $\pmb{\alpha}_K(u)$.   For any $K\in \mathcal{K}_o^n$, let  $\sigma_K\subseteq \partial K$ be the set given by
	\[
	\sigma_K=\big\{x\in\partial K: \   \pmb{\nu}_K(x) \ \ \mathrm{has\ more\ than\ one\ element}\big\}.
	\]  Denote by $\partial'K=\partial K\setminus\sigma_{K}$ the set of points on $\partial K$ that have a unique outer unit normal vector and by  $\mathcal{H}^{n-1}$  the $(n-1)$-dimensional Hausdorff  measure. According to \cite[p.84]{Sch},   $\mathcal{H}^{n-1}(\sigma_K)=0$ and hence $\mathcal{H}^{n-1}(\partial'K)=\mathcal{H}^{n-1}(\partial K)$. 
	
 The following lemma provides convenient formulas to calculate integrals with respect to the measure $\cV(K, \cdot).$ Recall that $\Psi_{K}(u)=\phi(\rho_K(u)u)[\rho_K(u)]^n$  for all $u\in\sphere.$
	\bl \label{two formula for dual orlicz curvature measure}
	Let $\phi$ be a function satisfying condition $C1)$.  For each $K\in\cK_o^n$, the following formulas
	\begin{eqnarray}
		\int_{S^{n-1}}g(v)d\cV(K,v)&=& \label{integral-general curvature measure 1}
		 \int_{S^{n-1}}g(\alpha_K(u))\Psi_{K}(u)du\\
		&=& \label{integral-general curvature measure 2}
		 \int_{\partial^\prime K}\langle x,\nu_K(x)\rangle g(\nu_K(x))\phi(x)\,d\mathcal{H}^{n-1}(x)
	\end{eqnarray}
	hold for any bounded Borel function $g: S^{n-1}\rightarrow \R$.
	\el

 \begin{proof}
		First, we prove  (\ref{integral-general curvature measure 1}). Let $\gamma(v)=\sum_{i=1}^m a_{i}\textbf{1}_{\eta_{i}}(v)$ for any $v\in \sphere$ be an arbitrary simple function, where $\eta_{i}\subseteq \sphere$ are Borel sets and   $\textbf{1}_A$ denotes the indicator function of the set $A$.  By \cite[(2.21)]{HLYZ}, one has  $u\in \pmb{\alpha}_K^*(\eta)$  if and only if $\alpha_K(u)\in \eta$, and this  further yields that  \begin{eqnarray*} \int_{S^{n-1}}  \gamma (\alpha_K(u))\Psi_{K}(u)\,du&=& \int_{S^{n-1}}\sum_{i=1}^m a_i \textbf{1}_{\eta_i}(\alpha_K(u))\Psi_{K}(u)\,du\\ &=& \int_{S^{n-1}}\sum_{i=1}^m a_i \textbf{1}_{\pmb{\alpha}^*_K(\eta_i)}(u)\Psi_{K}(u)\,du\\  &=& \sum_{i=1}^m a_i\int_{S^{n-1}} \textbf{1}_{\pmb{\alpha}^*_K(\eta_i)}(u)\Psi_{K}(u)\,du.
		\end{eqnarray*} Together with Definition \ref{general dual Orlicz curvature measure}, one has  \begin{eqnarray*} \int_{S^{n-1}}  \gamma (\alpha_K(u))\Psi_{K}(u)\,du
		&=& \sum_{i=1}^m a_i\int_{S^{n-1}} \textbf{1}_{\pmb{\alpha}^*_K(\eta_i)}(u)\Psi_{K}(u)\,du\\ &=&\sum_{i=1}^m a_i \cV(K,\eta_i)\\&=&\sum_{i=1}^m a_{i} \int_{S^{n-1}}\textbf{1}_{\eta_{i}}(v)\,d\cV(K,v)\\ &=&\int_{S^{n-1}}\gamma(v)\,d\cV(K,v).\end{eqnarray*} That is, (\ref{integral-general curvature measure 1}) holds true for simple functions. Following  from a standard limit approach by simple functions, one can prove formula (\ref{integral-general curvature measure 1}) for general bounded Borel functions $g: \sphere\rightarrow \R$.
		
		Next we  prove (\ref{integral-general curvature measure 2}). According to
		\cite[(2.31)]{HLYZ}, for each bounded integrable function $f:\sphere\rightarrow\R$, one has
		\begin{eqnarray*} \int_{S^{n-1}}f(u)\phi(\rho_K(u)u)\,du&=&\int_{\partial^\prime K}\langle x,\nu_K(x)\rangle f(\bar{x}) \frac{\phi(\rho_K(\bar{x})\bar{x})}{\rho_K^n(\bar{x})}\,d\mathcal{H}^{n-1}(x)\\ &=& \int_{\partial^\prime K}\langle  x,\nu_K(x)\rangle f(\bar{x})\frac{\phi(x)}{|x|^n} \,d\mathcal{H}^{n-1}(x), \end{eqnarray*}
		where $\bar{x}=x/|x|$, $\rho_K(\bar{x})\bar{x}=x$,  and $\rho_K(\bar{x})=|x|$.  Together with (\ref{integral-general curvature measure 1}) and the fact that $f=g\circ \alpha_K$ is bounded integrable on $S^{n-1}$, one has
		\begin{eqnarray*}
			\int_{S^{n-1}}g(v)\,d\cV(K,v)&=& \int_{S^{n-1}}g(\alpha_K(u))\Psi_{K}(u)\,du\\
			&=&  \int_{\partial^\prime K}\langle x,\nu_K(x)\rangle g(\nu_K(x))\phi(x)\,d\mathcal{H}^{n-1}(x).
		\end{eqnarray*} Hence, (\ref{integral-general curvature measure 2}) holds true.  \end{proof}
		
	The weak convergence of the general dual Orlicz curvature measure is proved in the following proposition.

	\bp \label{weak continuity} Let $\phi$  be a function satisfying condition $C1)$. If the sequence $\{K_i\}_{i=1}^{\infty} \subseteq \cK_o^n$ converges to $K\in \cK_o^n$ in the Hausdorff metric, then $\cV(K_i,\cdot)$ converges to $\cV(K,\cdot)$ weakly. \ep
		
	\begin{proof} As $\{K_i\}_{i=1}^{\infty}\subseteq  \cK_o^n$ converges to $K\in \cK_o^n$,   then $\rho_{K_i}$ converges to $\rho_{K}$ uniformly (see (\ref{conv-Hausdorff-radial})) and hence one can find constants $R_{1},R_{2}>0$, such that,  for all $u\in S^{n-1}$  and for  all $i\geq 1$, 
 \[ R_{1}\leq\rho_{K_{i}}(u)\leq R_{2} \ \ \mathrm{and} \ \ R_{1}\leq\rho_{K}(u)\leq R_{2}. \]
 For any fixed $u\in\sphere$ and for any function $\phi$ satisfying condition C1), it can be checked that  \begin{equation}\label{continuity for dual orlicz curvature measure 1} 
		 \Psi_{K_{i}}(u)=\phi(\rho_{K_i}(u)u)[\rho_{K_i}(u)]^n \rightarrow \phi(\rho_K(u)u) [\rho_K(u)]^n =\Psi_{K}(u) \ \  \mathrm{uniformly\ \ on}\ \ \sphere.
	\end{equation} 	
	
 Note that $\alpha_{K_i}\rightarrow \alpha_K$ almost everywhere on $S^{n-1}$ (see \cite[Lemma 2.2]{HLYZ}).  For any continuous function  $g: S^{n-1}\rightarrow \R$,
		by (\ref{continuity for dual orlicz curvature measure 1}), there exists a constant $M>0$, such that, for all $u\in \sphere$ and for all $i=1,2,\cdots$,  
		$$
		|g(\alpha_{K_i}(u))\Psi_{K_{i}}(u)|\leq M \ \ \ \ \mathrm{and}\ \ \ |g(\alpha_{K}(u))\Psi_{K}(u)| \leq M.
		$$ It follows from the dominated convergence theorem that
		\begin{eqnarray*}
			 \lim_{i\rightarrow \infty}\int_{S^{n-1}}g(\alpha_{K_i}(u))\Psi_{K_i}(u)\,du &=& \int_{S^{n-1}}\lim_{i\rightarrow \infty}g(\alpha_{K_i}(u))\Psi_{K_{i}}(u)\,du\\ &=& \int_{S^{n-1}}g(\alpha_{K}(u))\Psi_{K}(u)\,du.
		\end{eqnarray*} Together with (\ref{integral-general curvature measure 1}), then 
		\begin{eqnarray*}\lim_{i\rightarrow \infty}\int_{S^{n-1}}g(v)\,d\cV(K_i,v)&=&\lim_{i\rightarrow \infty}\int_{S^{n-1}}g(\alpha_{K_i}(u))\Psi_{K_i}(u)\,du\\  &=&\int_{S^{n-1}}g(v)\,d\cV(K,v),\end{eqnarray*} hold for any continuous function  $g: S^{n-1}\rightarrow \R$. In conclusion,  $\cV(K_i,\cdot)$ converges weakly to $\cV(K,\cdot)$ as desired. \end{proof}
	
Denote by $S_K(\cdot)$ the surface area measure of $K\in \cK_o^n$, namely, for any Borel set $\eta\subseteq\sphere,$
	\begin{equation}\label{surface:area:1}
	S_K(\eta)=\mathcal{H}^{n-1}(\nu_{K}^{-1}(\eta)),
	\end{equation}
	where $\nu_{K}^{-1}(\eta)$ is the reverse spherical image of $\eta$, i.e., $\nu_{K}^{-1}(\eta)=\{x\in \partial K: \ \ \nu_K(x)\in \eta\}.$ 
	\bp \label{ absolutely continuous}
	Let $K\in \cK_o^n$ and $\phi$ be a function satisfying condition C1). Then, $\cV(K,\cdot)$ is absolutely continuous with respect to the surface area measure $S_K(\cdot)$. \ep
	
	\begin{proof} Let $\eta\subseteq \sphere$ be a Borel set and $g=\pmb{1}_\eta$ in (\ref{integral-general curvature measure 2}). Then 
		\begin{eqnarray*}
			\cV(K,\eta)= \int_{\nu_K^{-1}(\eta)}\langle x,\nu_K(x)\rangle{\footnotesize }\phi(x)
			d\mathcal{H}^{n-1}(x).
		\end{eqnarray*}
		Since $K\in \cK_o^n$ and $\phi$ is a function satisfying condition C1), there exists a constant $T<\infty$, such that, $\langle x,\nu_K(x)\rangle \phi(x)\leq T$ for all $x\in\partial K$. Then
		\begin{eqnarray*}
			\int_{\nu_K^{-1}(\eta)}\langle x,\nu_K(x)\rangle \phi(x)
			d\mathcal{H}^{n-1}(x)\le T\int_{\nu_K^{-1}(\eta)}d\mathcal{H}^{n-1}(x).
		\end{eqnarray*}
		If $\eta \subseteq S^{n-1}$ is a Borel set such that $S_K(\eta)=0$, then $\mathcal{H}^{n-1}(\nu_K^{-1}(\eta))=0$ and thus
		\begin{eqnarray*}
			\cV(K,\eta)\le T\cdot \mathcal{H}^{n-1}(\nu_K^{-1}(\eta))=0.
		\end{eqnarray*}
		As a result, $\cV(K,\cdot)$ is absolutely continuous with respect to $S_K(\cdot)$.
	\end{proof}
	
Let us discuss the measure $\cV(K,\cdot)$ for $K\in \cK_o^n$ under Case 1 and Case 2 given in Section \ref{Section:2}. In Case 1, i.e., $\phi(x)=\psi(|x|)$, it follows from Definition \ref{general dual Orlicz curvature measure} that   for any Borel set $\eta\subseteq \sphere$, 
\begin{eqnarray}\cV(K, \eta)&=&\int_{\pmb{\alpha}^*_K(\eta)}\phi(\rho_K(u)u)[\rho_K(u)]^n\,du \nonumber \\ &=&\int_{\pmb{\alpha}^*_K(\eta)}\psi(\rho_K(u))[\rho_K(u)]^n\,du \nonumber \\ &=& \frac{1}{n}\int_{\pmb{\alpha}^*_K(\eta)}\hat{\varphi}(\rho_K(u))\,du,\label{special-curvature-measure--11}
\end{eqnarray} where $\hat{\varphi}(t)=n\psi(t)t^n$.  Recall that for $K\in \cK_o^n$ and $\varphi: (0, \infty)\rightarrow (0, \infty)$ a continuous function, the dual $L_{\varphi}$ Orlicz curvature measure of $K$, denoted by $\widetilde{C}_\varphi(K,\cdot)$, is defined in \cite{ZSY2017} as follows:  for each Borel set $\eta\subseteq \sphere$, 
	\begin{equation*}
		\widetilde{C}_\varphi(K,\eta)
		=\frac{1}{n}\int_{\pmb{\alpha}^*_K(\eta)}\varphi(\rho_K(u))du. 
	\end{equation*}  Hence, (\ref{special-curvature-measure--11}) asserts that $\cV(K, \cdot)=\widetilde{C}_{\hat{\varphi}}(K,\cdot)$.  In particular, if $$\phi(x)=\frac{|x|^{q-n}}{n}$$ which leads to $\hat{\varphi}(t)=t^q$, then $\cV(K, \cdot)$ is just the $q$-th dual curvature measure of $K$  \cite{HLYZ}; that is, for any Borel set $\eta\subseteq \sphere$, 
$$\cV(K, \eta)=\widetilde{C}_{\hat{\varphi}}(K,\eta)  =\frac{1}{n}\int_{\pmb{\alpha}^*_K(\eta)}[\rho_K(u)]^q \,du. $$
In Case 2,   i.e.,  $\phi(x)=\psi(|x|)\phi_{2}(\bar{x})$,  one has,   for any Borel set $\eta\subseteq \sphere$,
	\begin{eqnarray} 
	\widetilde{C}_{\phi,\mathcal{V}}(K,\eta)
	&=&\int_{\pmb{\alpha}^*_K(\eta)}\phi(\rho_K(u)u)[\rho_K(u)]^n\,du \nonumber \\&=& \int_{\pmb{\alpha}^*_K(\eta)} \psi(\rho_K(u))[\rho_K(u)]^n \phi_{2}(u)\,du \nonumber \\&=&\frac{1}{n}\int_{\pmb{\alpha}^*_K(\eta)}\hat{\varphi}(\rho_K(u)) \phi_{2}(u)\,du. \label{homogenes-2-2}
	\end{eqnarray}  In this case, Lemma \ref{two formula for dual orlicz curvature measure} can be rewritten as follows.
	\bc \label{measure-change}
	Let  $\phi(x)=\psi(|x|)\phi_{2}(\bar{x})$ satisfy condition C1). For $K\in\cK_o^n$, then
	\begin{eqnarray*}
		\int_{S^{n-1}}g(v)d\cV(K,v)&=&
		\frac{1}{n}\int_{S^{n-1}}g(\alpha_K(u))\hat\varphi(\rho_K(u))\phi_{2}(u)du\\
		&=&
		\frac{1}{n}\int_{\partial^\prime K}\langle x,\nu_K(x)\rangle \cdot g(\nu_K(x))\frac{\hat{\varphi}(|x|)\phi_{2}(\bar{x})}{|x|^n}\,d\mathcal{H}^{n-1}(x)\\&=&
		 \int_{\partial^\prime K}\langle x,\nu_K(x)\rangle \cdot g(\nu_K(x))\psi(|x|)\phi_{2}(\bar{x}) \,d\mathcal{H}^{n-1}(x),
	\end{eqnarray*}
	hold for each bounded Borel function $g: S^{n-1}\rightarrow \R$.
	\ec  
	\section{A variational interpretation for the general dual Orlicz curvature measure}
	
Let $f:\Omega\rightarrow (0, \infty)$ be a continuous function with $\Omega$  a closed set of $S^{n-1}$ such that $\Omega$ is not contained in any closed hemisphere of $S^{n-1}$. We shall need the Wulff shape and the convex hull of $f$,  denoted by $[f]$ and $\langle f\rangle$ respectively, whose definitions are given by  
	\begin{eqnarray*} 
	[f]= \cap_{u\in\Omega}
	\{x\in\R^n:\langle x, u\rangle\leq f(u)\} \ \ \mathrm{and}\ \ 
	\langle f\rangle =\text{conv}\big(\{f(u)u: u\in\Omega\}\big).
	\end{eqnarray*} Clearly, $[h_K]=K$ and $\langle \rho_K \rangle=K$ for each $K\in \cK_o^n$. 
 It can be easily checked that   $[f]\in\mathcal{K}_{o}^n$ and $\langle f\rangle\in\mathcal{K}_{o}^n$ for all  continuous functions $f:\Omega\rightarrow (0, \infty)$, due to the fact that $\Omega\subseteq \sphere$ is not contained in any closed hemisphere of $S^{n-1}$. As stated in \cite[Lemma 2.8]{HLYZ},
	\begin{equation}\label{relation}
		[f]^*=\langle 1/f \rangle.
	\end{equation} 	When $\Omega=S^{n-1}$ and $f\in C^+(\sphere)$, the set of all positive continuous functions defined on $\sphere$, it is obvious that $h_{[f]}(u)\le f(u)$ for all $u\in \sphere$. A less trivial fact is that $h_{[f]}(u)= f(u)$ for almost all $u\in\sphere$ with respect to the surface area measure $S_{[f]}(\cdot)$.  	
	
The variational interpretation for the general dual Orlicz curvature measure is stated as follows. Let $\Omega$ be a closed set of $S^{n-1}$ such that $\Omega$ is not contained in any closed hemisphere of $S^{n-1}$. \bt\label{variational-for-decreaing-1} Let $h_0: \Omega\rightarrow (0,\infty)$ and $g: \Omega\rightarrow \R$ be two continuous functions. Define $h_{t}$ by 
		\begin{equation}\label{the definition of ht}
			\log(h_{t}(u))=\log(h_{0}(u))+tg(u)+o(t,u)  \ \  \mathrm{for\  all}\  u\in\Omega,
		\end{equation}
		where $o(t,\cdot):\Omega\rightarrow\R$ is continuous  and  $o(t,u)/t \rightarrow 0$ uniformly on $\Omega$ as $t\rightarrow0.$  Let $\phi$ be a function satisfying condition C1). Then  \begin{eqnarray}\label{variational formula-2}\frac{d}{dt}\eV([h_t])\bigg |_{t=0} = -\int_\Omega g(u)\,d\cV([h_0],u).
	\end{eqnarray} \et
	
\vskip 2mm \noindent {\bf Remark.} An immediate consequence of (\ref{variational formula-2}) and the chain rule for derivative is the following formula, which will be used in solving the general dual Orlicz-Minkowski problem:   \begin{eqnarray} \label{variational formula-3}
		\frac{d}{dt}\log\eV([h_t])\bigg |_{t=0} = - \frac{1}{\eV([h_0])}\int_\Omega g(u)\, d\cV([h_0],u).\end{eqnarray}

	\begin{proof} Let $\rho_{0}:\Omega\rightarrow (0, \infty) $ be a continuous function. For $\delta>0$ and $t\in(-\delta,\delta)$, let
		\begin{equation*} 
			\log(\rho_{t}(u))=\log(\rho_{0}(u))+tg(u)+o(t,u) \ \ \mathrm{for\ all}\ u\in\Omega,
		\end{equation*}
		where $o(t,\cdot):\Omega\rightarrow\R$ is  continuous and $o(t,u)/t \rightarrow 0$ uniformly on $\Omega$ as $t\rightarrow0.$

		First of all,  let us prove the following formula: for almost every $u\in\sphere$ (with respect to the spherical measure), 
		\begin{eqnarray}\label{variation formula bounded 1}
			\frac{d}{dt}\Phi(\rho_{\langle\rho_{t}\rangle^*}(u),u)\Big|_{t=0} = \frac{d}{dt} \int_{\rho_{\langle\rho_{t}\rangle^*}(u)}^{\infty}\phi(ru)r^{n-1}\,dr \Big|_{t=0} = \Psi_{\langle\rho_{0}\rangle^{*}}(u)	
			g(\alpha^*_{\langle \rho_0\rangle}(u)).
		\end{eqnarray} In fact, it follows from the chain rule  and $\rho_{\langle\rho_{t}\rangle^*}(u)=h_{\langle\rho_{t}\rangle}^{-1}(u)$ for all $u\in \sphere$ that
		\begin{eqnarray*}
			\frac{d}{dt}\Phi(\rho_{\langle\rho_{t}\rangle^*}(u),u)\Big |_{t=0} &=&\frac{d}{dt}\int_{e^{-\log h_{\langle\rho_{t}\rangle}(u)}}^{\infty}\phi(ru)r^{n-1}\,dr  \Big |_{t=0}\\
			&=&\phi(h_{\langle\rho_{0}\rangle}^{-1}(u)u)h_{\langle\rho_{0}\rangle}^{-n}(u)\cdot \frac{d}{dt}\log h_{\langle\rho_{t}\rangle}(u)\Big|_{t=0}\\
			&=&\Psi_{\langle\rho_{0}\rangle^{*}}(u) \cdot	
			g(\alpha^*_{\langle \rho_0\rangle}(u)),
		\end{eqnarray*}  where the last equality follows from  \cite[(4.4)]{HLYZ}, i.e.,   
		\begin{equation*} 
			\lim_{t\rightarrow 0}\frac{\log h_{\langle \rho_t\rangle}(v)-\log h_{\langle \rho_0\rangle}(v)}{t}=g(\alpha^*_{\langle \rho_0\rangle}(v))\end{equation*} holds  for any $v\in S^{n-1}\setminus \eta_0$, with  $\eta_0 =\eta_{\langle \rho_0\rangle}$ the complement of the set  of the regular normal vectors of $\langle \rho_0\rangle$. Note that the spherical measure of $\eta_0$ is zero.  
			
We shall need the following argument in order to use the dominated convergence theorem:   there exist two constants $\delta>0$ and $M>0$, such that, for all $t\in (-\delta, \delta)$ and for all $u\in S^{n-1}$,
		\begin{eqnarray}\label{bounded-1-2-1}
			 \big| \Phi(\rho_{\langle\rho_{t}\rangle^*}(u),u)-\Phi(\rho_{\langle\rho_{0}\rangle^*}(u),u) \big|\le M|t|.
		\end{eqnarray} Note that $\langle\rho_{t}\rangle\rightarrow \langle\rho_{0}\rangle$ in the Hausdorff metric; this is a direct consequence of the Aleksandrov's convergence lemma \cite{Aleks1938} and formula (\ref{conv-Hausdorff-radial}). Therefore, $\rho_{\langle\rho_{t}\rangle^{*}}\rightarrow \rho_{\langle\rho_{0}\rangle^{*}}$ uniformly on $\sphere$. As $\langle\rho_{0}\rangle^{*}\in\cK_o^n$, one can find constants $l_{1}, l_{2}, \delta_{1}>0$, such that, $
l_{1}<\rho_{\langle\rho_{t}\rangle^{*}}(u)<l_{2}$ holds for all $u\in\sphere$ and for all $t\in(-\delta_{1}, \delta_{1}).$  It follows from condition C1) and the continuity of $\phi$ that 
	\begin{eqnarray} \big|\big[\log \Phi(e^{-s},u)\big]'\big|=\big|\phi(e^{-s}u)e^{-sn}/\Phi(e^{-s},u)\big|\leq L_2 \label{uniform-bound--1} \end{eqnarray}
		holds for some finite constant $L_2$ independent of $u\in\sphere$ and for all $s\in (-\log l_2, -\log l_1)$. Note that $\log h_{\langle \rho_t\rangle}(u)\in (-\log l_2, -\log l_1)$ and $\log h_{\langle \rho_0\rangle}(u)\in (-\log l_2, -\log l_1)$ for all $u\in \sphere$ and for all $t\in(-\delta_{1}, \delta_{1}).$ By (\ref{uniform-bound--1}) and the mean value theorem,  one has, 
		for all $u\in\sphere$ and for all $t\in(-\delta,\delta)$ (without loss of generality, we can assume that  $0<\delta<\delta_1$),  \begin{eqnarray} \label{bounded-2-2-1-3}  \Big|\log \Phi(h_{\langle\rho_{t}\rangle}^{-1}(u),u)-\log\Phi(h_{\langle\rho_{0}\rangle}^{-1}(u),u)\Big| \leq L_2 \ \Big|\log h_{\langle \rho_t\rangle}(u)-\log h_{\langle \rho_0\rangle} (u)\Big|\le L_{2}M_{1}|t|,
		\end{eqnarray} where the last inequality follows from \cite[Lemma 4.1]{HLYZ}, i.e, there exist constants $0<\delta, M_1<\infty$ such that, for all $u\in\sphere$ and for all $t\in(-\delta,\delta)$, \begin{equation*} 
			\Big|\log h_{\langle \rho_t\rangle}(u)-\log h_{\langle \rho_0\rangle}(u)\Big |\le M_{1}|t|.
		\end{equation*}	

It follows from condition C1) that there is a constant $L_1$ (independent of $u\in \sphere$), such that, for all $u\in \sphere$ and for all $t\in(-\delta,\delta)$,
		\[ 0<s_{t}=\frac{\Phi(\rho_{\langle\rho_{t}\rangle^{*}}(u),u)}{\Phi(\rho_{\langle\rho_{0}\rangle^{*}}(u),u)}=\frac{\Phi(h_{\langle\rho_{t}\rangle}^{-1}(u),u)}{\Phi(h_{\langle\rho_{0}\rangle}^{-1}(u),u)}<L_1.\] Hence $|s_t-1|\le L_1 \cdot |\log s_t|$ (see e.g., \cite[p.362]{HLYZ}). Together with inequality (\ref{bounded-2-2-1-3}), one gets, for all $u\in\sphere$ and for all $t\in(-\delta,\delta)$,
		\begin{eqnarray*} \Big| \Phi(\rho_{\langle\rho_{t}\rangle^*}(u),u)-\Phi(\rho_{\langle\rho_{0}\rangle^*}(u),u) \Big| &=&
			\Big|\Phi(h_{\langle\rho_{t}\rangle}^{-1}(u),u)-\Phi(h_{\langle\rho_{0}\rangle}^{-1}(u),u)\Big|\\
			&\leq& \Phi(h_{\langle\rho_{0}\rangle}^{-1}(u),u)\cdot L_1 \cdot \Big|\log \Phi(h_{\langle\rho_{t}\rangle}^{-1}(u),u)-\log\Phi(h_{\langle\rho_{0}\rangle}^{-1}(u),u)\Big| \\ &\leq & \Phi(h_{\langle\rho_{0}\rangle}^{-1}(u),u)\cdot L_1  L_2 M_1\cdot |t| \\ &\leq & \Phi(l_1,u)\cdot L_1  L_2 M_1\cdot |t|. \end{eqnarray*} That is, inequality (\ref{bounded-1-2-1}) holds by letting  $M=   L_1  L_2 M_1\cdot \max_{u\in\sphere} \Phi(l_1, u)<\infty.$ 
		
Now we are ready to prove formula (\ref{variational formula-2}). To this end, let  $[h_t]$ be the Wulff shape associated to $h_t$ with $h_t$ given by (\ref{the definition of ht}). Consider   $\kappa_t=1/h_t$ and then  $$\log \kappa_t=-\log h_t=-\log h_{0}-tg-o(t,\cdot)=\log \kappa_0-t g- o(t,\cdot).$$ Moreover, $[h_t]=\langle 1/h_t\rangle^*=\langle \kappa_t\rangle^*$ due to the bipolar theorem and (\ref{relation}).  It follows from (\ref{general-quermass-33}), (\ref{variation formula bounded 1}),  (\ref{bounded-1-2-1}) and the dominated convergence theorem that \begin{eqnarray*}
			\frac{d}{dt}\eV([h_t])\Big|_{t=0}&=&\frac{d}{dt}\eV(\langle\kappa_t\rangle^*)\Big|_{t=0}\\
			&=& \frac{d}{dt}\int_{S^{n-1}}\Phi(\rho_{\langle\kappa_t\rangle^*}(u),u)\, du\Big|_{t=0}\\&=& \int_{S^{n-1}}\frac{d}{dt}\Phi(\rho_{\langle\kappa_t\rangle^*}(u),u)\Big |_{t=0}\, du\\&=&- \int_{S^{n-1}}\Psi_{\langle\kappa_{0}\rangle^{*}}(u)	\cdot 
			g(\alpha^*_{\langle \kappa_{0}\rangle}(u))\, du.
		\end{eqnarray*}  Together with  (\ref{integral-general curvature measure 1}) and the fact that the spherical measure of $\eta_{0}$ is zero, one can prove formula (\ref{variational formula-2}) as follows:  \begin{eqnarray*}
			\frac{d}{dt}\eV([h_t])\Big|_{t=0}&=&- \int_{S^{n-1}\setminus \eta_0} \Psi_{\langle\kappa_{0}\rangle^{*}}(u)\cdot 	
			g(\alpha^*_{\langle \kappa_{0}\rangle}(u))\, du\\&=&- \int_{S^{n-1}} (\hat{g}\textbf{1}_\Omega)(\alpha_{\langle \kappa_{0}\rangle^*}(u))\Psi_{\langle\kappa_{0}\rangle^{*}}(u)\, du\\ \nonumber
			&=&- \int_{S^{n-1}} (\hat{g}\textbf{1}_\Omega)(u)\, d\cV(\langle \kappa_{0}\rangle^*,u)\\&=&-\int_\Omega g(u)\,d\cV([h_0],u),
		\end{eqnarray*}
		where $\hat{g}: S^{n-1}\rightarrow \R$ is a continuous function,  such that, for all $v\in S^{n-1}\setminus \eta_0$, \begin{eqnarray*}\label{extended function}
			g(\alpha_{\langle \rho_0\rangle^*}(v))=(\hat{g}\textbf{1}_\Omega)(\alpha_{\langle \rho_0\rangle^*}(v)).
		\end{eqnarray*}
		The existence of such $\hat{g}$ was proved in  \cite[p.364]{HLYZ}.\end{proof} 
		
	\section{A solution of the general dual Orlicz-Minkowski problem}
	
	In this section, we provide a solution to the following general dual Orlicz-Minkowski problem. 
	
	\vskip 2mm \noindent {\bf The general dual Orlicz-Minkowski problem:} 
	{\em Given a nonzero finite Borel measure $\mu$ defined on $\sphere$ and a continuous function $\phi: \Rn\setminus\{o\}\rightarrow (0, \infty)$,  can one find a constant  $\tau>0$ and a convex body $K$ (ideally $K\in \cK_{o}^n$), such that,  $\mu=\tau\cV(K,\cdot)?$} 
	 
	 Clearly, if the  general dual Orlicz-Minkowski problem has solutions, the constant $\tau$ can be calculated by $$|\mu|=\int_{\sphere} \,d\mu(v)=\tau \int_{\sphere} \,d\cV(K, v)=\tau\cdot \cV(K, \sphere)$$ and equivalently \begin{equation}\label{constant--tau--1} \tau=\frac{|\mu|}{\cV(K, \sphere)}. \end{equation}

 It is well known that, to have the various Minkowski problems solvable, the given measure $\mu$ must satisfy  that  $\mu$ is not concentrated in any closed hemisphere, i.e.,
\begin{equation}\label{condition for Minkowski problem}
		\int _{\sphere} \langle \xi, \theta\rangle_+ \,d\mu(\theta)>0\ \ \ \mathrm{for\ \  all\ \ } \xi\in \sphere.
	\end{equation} Hereafter, $a_+=\max\{a, 0\}$ for $a\in \R$. 
	
	In fact, (\ref{condition for Minkowski problem}) is also a necessary condition in our setting. That is, if there exists a convex body $K\in \cK_{o}^n$, such that, $$\frac{\mu}{|\mu|}=\frac{\widetilde{C}_{\phi,\mathcal{V}}(K,\cdot)}{\widetilde{C}_{\phi,\mathcal{V}}(K,S^{n-1})}, $$
	then $\mu$ satisfies (\ref{condition for Minkowski problem}). 
	To this end, let $\xi\in\sphere$ be given. Then \begin{equation}
		\int _{\sphere} \langle \xi, v\rangle_+ \,d\mu(v)=\frac{|\mu|}{\cV(K,\sphere)}\int_{\sphere}\langle \xi, v\rangle_+\,d\cV(K,v).\label{concentration-0--10}
	\end{equation}
	Hence, in order to show that $\mu$ satisfies (\ref{condition for Minkowski problem}), it is enough to show that
	\[
	\int_{\sphere}\langle \xi, v\rangle_+\,d\cV(K,v)>0.
	\]
	In fact, it follows from (\ref{integral-general curvature measure 2}) that 
	\begin{eqnarray*}
		\int_{S^{n-1}}\langle \xi, v\rangle_+\,d\cV(K, v)
		&=&\int_{\partial^\prime K}  \langle\xi, \nu_K(x)\rangle_+\cdot \langle x,\nu_K(x)\rangle \phi(x)\,d\mathcal{H}^{n-1}(x).
	\end{eqnarray*}	
	As $K\in\mathcal{K}_{o}^n$, one can find a constant $M$ such that $\langle x, \nu_{K}(x)\rangle\phi(x)\geq M$ for all $x\in\partial' K$. Consequently,
	\begin{eqnarray}
		\int_{\sphere}\langle \xi, v\rangle_+\,d\cV(K, v)
		\geq M\int_{\partial 'K}\langle \xi,\nu_{K}(x)\rangle_{+}\,d\mathcal{H}^{n-1}(x)>0,\label{concentration-1--1}
	\end{eqnarray}	
	as the surface area measure $S_K(\cdot)$ satisfies $$\int_{\partial' K}\langle \xi,\nu_{K}(x)\rangle_{+}\,d\mathcal{H}^{n-1}(x)=\int_{\partial K}\langle \xi,\nu_{K}(x)\rangle_{+}\,d\mathcal{H}^{n-1}(x)=\int_{\sphere}\langle \xi,u\rangle_{+}\,dS_K(u)>0.$$
		
 The following theorem also shows that (\ref{condition for Minkowski problem}) is a sufficient condition for the general dual Orlicz-Minkowski problem.

	\bt \label{solution-general-dual-Orlicz-main theorem}  Let $\mu$ be a nonzero finite Borel measure on $S^{n-1}$ satisfying (\ref{condition for Minkowski problem}) and let $\phi$ be a function satisfying conditions C1) and C2).
	Then there exists a convex body $K\in \cK_o^n$, such that, $$\frac{\mu}{|\mu|}=\frac{\widetilde{C}_{\phi,\mathcal{V}}(K,\cdot)}{\widetilde{C}_{\phi,\mathcal{V}}(K,S^{n-1})}. $$
	\et
	
	In order to prove Theorem \ref{solution-general-dual-Orlicz-main theorem}, we need the following lemma. \bl\label{bounded-lemma-Minkowski-2} Let $\mu$ be a nonzero finite Borel measure on $S^{n-1}$ satisfying (\ref{condition for Minkowski problem}) and let $\phi$ be a function satisfying conditions C1) and C2). Then there exists a convex body $Q_0\in\cK_o^n$ such that $\eV(Q_0)=|\mu|$ and
	\begin{equation}\label{optimization problem}
		\mathcal{F}(Q_0)=\sup\big\{\mathcal{F}(K): \eV(K)=|\mu| \ \mbox{and}\ K\in \cK_o^n \big\},
	\end{equation}  where  $\mathcal{F}: \cK_o^n\rightarrow \R$ is defined  by
	\begin{equation}\label{function-optimal-1}
		\mathcal{F}(K)=-\frac{1}{|\mu|}\int_{S^{n-1}}\log h_K(v)d\mu(v).
	\end{equation}
	\el
	\begin{proof} Let $\{Q_i\}_{i=1}^{\infty}\subseteq \cK_o^n$ be  such that $\eV(Q_i)=|\mu|$ and
		\begin{equation}\label{maximal-seq-1}
		\lim_{i\rightarrow \infty} \mathcal{F}(Q_i)=\sup\Big\{\mathcal{F}(K): \eV(K)=|\mu| \ \text{and}\ K\in \cK_o^n\Big\}.\end{equation}   First of all, we claim that the sequence $\{Q_i^{*}\}_{i=1}^\infty$ is uniformly bounded.  That is, we need to prove that there exists a constant $R>0$ such that $
		Q_i^{*} \subseteq R \ball$ for all $i=1,2,\cdots$
 
	Assume not, i.e., there are no finite constants $R$ such that  $Q_i^{*} \subseteq R \ball$  for all $i=1,2,\cdots$.	 Let $v_i\in \sphere$ be such that $\rho_{Q_i^{*}}(v_{i})=\max_{u\in\sphere}\rho_{Q_{i}^{*}}(u)$ and $R_{Q_i^{*}}=\rho_{Q_i^{*}}(v_{i})$. Without loss of generality,  we can assume that $R_{Q_i^{*}}\rightarrow\infty$  (otherwise, the sequence $\{Q_i^{*}\}_{i=1}^\infty$ is uniformly bounded) and $v_{i}\rightarrow v_{0}$ (due to the compactness of $\sphere$) as $i\rightarrow\infty$. Consequently, for any  $M>0,$ there exists $i_{M}>0$ such that $R_{Q_i^{*}}\geq M$ for all $i>i_{M}.$ Clearly,  for all $i>i_{M},$
		\begin{equation*}
			h_{Q_i^{*}}(u)\ge \langle u, v_i\rangle_+\ R_{Q_i^{*}}\geq M \langle u, v_i\rangle_+. 
		\end{equation*}	 Recall that $\Phi(t,u)=\int_{t}^{\infty}\phi(ru)r^{n-1}dr$ is decreasing on $t$. Then for all $i>i_{M}$ and for all $u\in\sphere,$
		\begin{equation}\label{the bound for lemma of solution}
			\Phi(\rho_{Q_{i}}(u),u)=\Phi(h_{Q_{i}^{*}}^{-1}(u),u)\geq\Phi([M\langle u,v_{i}\rangle_{+}]^{-1},u),
		\end{equation}	
		where we let $\Phi([M\langle u,v_{i}\rangle_{+}]^{-1},u)=0$ if $\langle u,v_{i}\rangle_{+}=0.$ 
		
		Fatou's lemma implies that \begin{eqnarray*}\nonumber 
			\liminf_{i\rightarrow \infty} \int_{S^{n-1}}\Phi([M\langle u,v_{i}\rangle_{+}]^{-1},u)\, du
			&=&\liminf_{i\rightarrow \infty} \int_{S^{n-1}}\int_{[M\langle u,v_{i}\rangle_{+}]^{-1}}^{\infty}\phi(ru)r^{n-1}\,dr\, du\\  \nonumber
			&\geq& \int_{S^{n-1}}\liminf_{i\rightarrow \infty}\int_{0}^{\infty}\pmb{1}_{([M\langle u,v_{i}\rangle_{+}]^{-1},\infty)}\phi(ru)r^{n-1}\,dr\, du\\ \nonumber
			&\geq& \int_{\sphere}\int_{0}^{\infty}\liminf_{i\rightarrow \infty}\pmb{1}_{([M\langle u,v_{i}\rangle_{+}]^{-1},\infty)}\phi(ru)r^{n-1}\,dr\, du\\ \nonumber
			&=& \int_{\sphere}\int_{[M\langle u,v_{0}\rangle_{+}]^{-1}}^{\infty}\phi(ru)r^{n-1}\,dr\, du\\
			&=& \int_{S^{n-1}}\Phi([M\langle u,v_{0}\rangle_{+}]^{-1},u)\, du.
		\end{eqnarray*}	Together with (\ref{general-quermass-33}) and (\ref{the bound for lemma of solution}), one has  
		\begin{eqnarray}\nonumber |\mu|&=&
			\lim_{i\rightarrow \infty}\eV(Q_i)\nonumber \\ &=&\lim_{i\rightarrow \infty} \int_{S^{n-1}}\Phi(\rho_{Q_{i}}(u),u)\, du \nonumber \\ &\geq&\liminf_{i\rightarrow \infty} \int_{S^{n-1}}\Phi([M\langle u,v_{i}\rangle_{+}]^{-1},u)\, du \nonumber\\ &\geq& \int_{S^{n-1}}\Phi([M\langle u,v_{0}\rangle_{+}]^{-1},u)\, du. \label{new-bound--1}\end{eqnarray}

For all  $j\geq 2$,  let $$\Sigma_j(v_0):=\Big\{u\in S^{n-1}:\langle u, v_0\rangle_+ >1/j\Big\}.$$ It follows from the monotone convergence theorem and the fact $\Sigma_j(v_0)\subseteq \Sigma_{j+1}(v_0)\subseteq \cup_{j=1}^{\infty}\Sigma_j(v_0)=\sphere\setminus\{u\in \sphere: \ \langle u, v_0\rangle=0\}$ that  
		\begin{eqnarray*} \lim_{j\rightarrow \infty}\int_{\Sigma_j(v_0)}  \langle u, v_0\rangle_+ \,du =\int_{\cup_{j=1}^{\infty}\Sigma_j(v_0)} \langle u, v_0\rangle_+ \,du = \int_{S^{n-1}} \langle u, v_0\rangle_+ \,du>0,
		\end{eqnarray*} where the last inequality is due to the fact that the spherical measure is not concentrated on any closed hemisphere. Hence, there exists $j_0\geq 2$, such that, \begin{equation*}
			 \int_{\Sigma_{j_0}(v_0)}\,du\geq \int_{\Sigma_{j_0}(v_0)}  \langle u, v_0\rangle_+\,
			du\geq\frac{1}{2}\int_{\sphere}\langle u, v_0\rangle_+\,
			du>0. \end{equation*} It can be checked that    $[M\langle u, v_0\rangle_{+}]^{-1}\leq j_0/M$  for all $u\in \Sigma_{j_0}(v_0)$. By (\ref{general-quermass-22}) and (\ref{new-bound--1}), one gets  
		\begin{eqnarray*}\nonumber |\mu| \geq
			 \int_{S^{n-1}}\Phi([M\langle u,v_{0}\rangle_{+}]^{-1},u)\, du \nonumber
			\geq \int_{\Sigma_{j_0}(v_0)}\Phi(j_0/M,u)\, du=\eV\big(\mathbb{R}^n\setminus \mathcal{C}(v_{0}, j_{0}/M,1/j_{0})\big), \nonumber
		\end{eqnarray*} where for any fixed $u_{0} \in \sphere$,   $a>0$ and $b_{0}\in (0, 1)$,  $$\mathcal{C}(u_{0},a,b_{0})=\Big\{x\in\mathbb{R}^n:\ \langle \bar{x},  u_{0}\rangle\geq b_{0}\ \mathrm{and}\ |x|\geq a\Big\}.$$ As $\phi$ satisfies condition C2),  one gets a contradiction as follows: $$\infty> |\mu| \geq \lim_{M\rightarrow\infty} \eV(\mathbb{R}^n\setminus \mathcal{C}\big(v_{0}, j_{0}/M,1/j_{0})\big) =\infty.$$  Therefore, the sequence $\{Q_i^{*}\}_{i=1}^{\infty}$  is uniformly bounded.

Without loss of generality, we assume that $Q_i^{*}\rightarrow Q$ (more precisely, a subsequence of $\{Q_i^{*}\}_{i=1}^{\infty}$) in the Hausdorff metric for some compact convex set $Q\subseteq \R^n$, due to the Blaschke selection theorem (see e.g., \cite{Sch}). Note that $Q$ may not be a convex body, however,  the support function of $Q$ can be defined as in  (\ref{support-function--1}) and  $Q_i^{*}\rightarrow Q$ in the Hausdorff metric is defined as in (\ref{H-metric-1}). 

We now  show $Q\in\mathcal{K}_{o}^n$ and the proof can be obtained by an argument almost identical to those in \cite{zhao, ZSY2017}. In fact, assume that $Q\notin \cK_o^n$ and $o\in \partial Q$.   Then, there exists $u_0\in S^{n-1}$ such that $ 
		\lim_{i\rightarrow \infty} h_{Q_i^{*}}(u_0)=h_{Q}(u_0)=0.$ Let $$
		\Sigma_{\delta_0} (u_0)=\{v\in S^{n-1}: \langle v,u_0\rangle >\delta_0\}.$$ By  (\ref{function-optimal-1}),  $\eV(Q_i)=|\mu|$ and $Q_i^{*} \subseteq R \ball$ (without loss of generality, let $R>1$) for all $i$,  one has
		\begin{eqnarray*}
			\mathcal{F}(Q_i)&=&-\frac{1}{|\mu|}\int_{S^{n-1}}\log h_{Q_i}(v)\,d\mu(v)\\ &=&\frac{1}{|\mu|}\int_{\Sigma_{\delta_0}(u_0)}\log \rho_{Q_i^*}(v)\,d\mu(v)+\frac{1}{|\mu|}\int_{S^{n-1}\setminus \Sigma_{\delta_0}(u_0)}\log \rho_{Q_i^*}(v)\,d\mu(v)\\ &\leq& \frac{1}{|\mu|}\int_{\Sigma_{\delta_0}(u_0)}\log \rho_{Q_i^{*}}(v)\,d\mu(v)+\log R.
		\end{eqnarray*} It follows from $\mu(\Sigma_{\delta_0}(u_0))>0 $ and  $\rho_{Q_i^{*}}\rightarrow 0$ on $\Sigma_{\delta_0} (u_0)$ uniformly for some $\delta_0>0$ that   $$\lim_{i\rightarrow \infty} \mathcal{F}(Q_i) =-\infty,$$  which is impossible. Hence, $o\in \text{int} Q$ and  then $Q\in \cK_o^n$.

	Finally, let us check that $Q_0=Q^{*}\in \cK_o^n$ satisfies $\eV(Q_0)=|\mu|$ and (\ref{optimization problem}). In fact, as $Q_i^{*}\rightarrow Q$, one has $Q_i\rightarrow Q^*=Q_0$ due to the bipolar theorem. Then  $$\eV(Q_0)=\lim_{i\rightarrow \infty} \eV(Q_i)=|\mu|$$ is an immediate consequence of Lemma \ref{continuity-general-dual-qu-1}.  On the other hand, $h_{Q_i} \rightarrow h_{Q_0}$  uniformly on $\sphere$  due to $Q_i\rightarrow Q_0\in \cK_o^n$ and  (\ref{H-metric-1}). Moreover, there exist constants $R_1, R_2\in (0, \infty)$, such that,  for all $u\in \sphere$ and for all $i\geq 1$, $$R_1\leq h_{Q_i}(u)\leq R_2 \ \  \mathrm{and} \ \  R_1\leq h_{Q_0}(u)\leq R_2.$$ These further imply that, for all $u\in \sphere$ and for all $i\geq 1$, $$|\log h_{Q_i}(u)|\leq \max\{ |\log R_1|, |\log R_2|\}<\infty.$$ It follows from the dominated convergence theorem that	\begin{eqnarray*}
		\lim_{i\rightarrow \infty} \mathcal{F}(Q_i)&=&\lim_{i\rightarrow \infty} -\frac{1}{|\mu|}\int_{S^{n-1}}\log h_{Q_i}(v)\,d\mu(v)\\ &=& -\frac{1}{|\mu|}\int_{S^{n-1}} \lim_{i\rightarrow \infty} \log h_{Q_i}(v)\,d\mu(v)\\ &=&  -\frac{1}{|\mu|}\int_{S^{n-1}}\log h_{Q_0}(v)\,d\mu(v)\\ &=&  \mathcal{F}(Q_0).
	\end{eqnarray*} Together with (\ref{maximal-seq-1}), one can easily get  the desired formula  (\ref{optimization problem}).  \end{proof}
\vskip 2mm \noindent {\it Proof of Theorem \ref{solution-general-dual-Orlicz-main theorem}.} Recall that each $K\in\mathcal{K}_{o}^n$ can be uniquely determined by its support function and vice versa. Thus we can let $\eV(h_{[f]})=\eV([f])$ for all $f\in C^+(S^{n-1})$. On the other hand, as $f\geq h_{[f]}$ for all $f\in C^+(S^{n-1})$, then
	 \begin{equation}\label{functional-functional-1}
	\mathcal{F}(f) :=-\frac{1}{|\mu|}\int_{S^{n-1}}\log f(v)\,d\mu(v)\le \mathcal{F}(h_{[f]}).
		\end{equation}   Consider the following optimization problem: \begin{equation}\label{max-problem-2}
		\sup\Big\{\mathcal{F}(f): \eV([f])=|\mu| \ \text{for}\ f\in C^+(S^{n-1})\Big\}.
	\end{equation} According to (\ref{functional-functional-1}) and  Lemma \ref{bounded-lemma-Minkowski-2}, the support function of convex body $Q_0\in \cK_o^n$ found in Lemma \ref{bounded-lemma-Minkowski-2} is an optimizer for the optimization problem (\ref{max-problem-2}).  
	
	On the other hand, the method of Lagrange multipliers can be used to find the necessary conditions for the optimizers for the optimization problem (\ref{max-problem-2}).  In fact, for $\delta>0$ small enough, let $h_t(v)=h_{Q_0}(v) e^{tg(v)}$ for $t\in (-\delta, \delta)$ and for $v\in S^{n-1}$, where $g: \sphere\rightarrow \R$ is an arbitrary continuous function.  Let  $$\mathcal{L}(t, \tau)=\mathcal{F}(h_t)-\tau \big(\log \eV([h_t]) -\log |\mu|\big).$$ As  $h_{Q_0}$ is an optimizer to (\ref{max-problem-2}), the following equation holds:
	\begin{equation}\label{lagrange method}
		\frac{\partial}{\partial t}\mathcal{L}(t, \tau)\Big|_{t=0}=0.
	\end{equation} It is easily checked that 
	\begin{eqnarray*}
		\frac{\partial}{\partial t}\mathcal{F}(h_t)\Big|_{t=0}=\frac{\partial}{\partial t}\bigg(\!-\frac{1}{|\mu|}\int_{S^{n-1}}[\log h_{Q_0}(v)+t g(v)]\,d\mu(v)\bigg)\bigg|_{t=0}
		=-\frac{1}{|\mu|}\int_{S^{n-1}}g(v)\,d\mu(v).
	\end{eqnarray*}	 It follows from (\ref{variational formula-3}) that 
	\begin{eqnarray*}
		\frac{\partial}{\partial t}\log \eV([h_t]) \Big|_{t=0}=-\frac{1}{\eV(Q_0)}\int_{S^{n-1}}g(v)\,d\cV(Q_0, v).
	\end{eqnarray*} Due to $\eV(Q_0)=|\mu|$, one can rewrite (\ref{lagrange method}) as follows:
	\begin{equation*} 
		\int_{S^{n-1}}g(v)\,d\mu(v)
		= \tau \int_{S^{n-1}}g(v)\,d\cV(Q_0, v)
	\end{equation*} holding for arbitrary continuous function $g: \sphere\rightarrow \R$. Consequently,  $\mu=\tau \cV(Q_0,\cdot)$  with the constant $\tau$  given by (\ref{constant--tau--1}), that is,   $$\tau=\frac{|\mu|}{\cV(Q_0,S^{n-1})}.$$ In summary, a solution to the general dual Orlicz-Minkowski problem has been found.     \hfill  $\Box$
	
	The following corollary provides a solution to the general dual Orlicz-Minkowski problem under the Case 2 in Section \ref{Section:2}, i.e.,   $\phi(x)=\psi(|x|)\phi_{2}(\bar{x})$ with $\psi: (0, \infty)\rightarrow (0, \infty)$ and $\phi_{2}: \sphere \rightarrow (0, \infty)$ continuous functions. Again, let $\hat{\phi}$ and $\psi$ be given as in (\ref{definition:psi:hat}) or (\ref{def:hat:psi-phi}), and $\hat{\varphi}(t)=  n\psi(t)t^n.$    
	\bc \label{cor-5-1} Let $\phi(x)=\psi(|x|)\phi_{2}(\bar{x})$ be a continuous function such that the continuous function $\phi_2$ is positive on $\sphere$, and the functions $\hat{\phi}$ and $\hat{\varphi}$ satisfy conditions A1)-A3). Then the following are equivalent: 
	\begin{itemize} \item[i)]  $\mu$ is a nonzero finite Borel measure on $S^{n-1}$ satisfying (\ref{condition for Minkowski problem});
	
\item [ii)] there exists a convex body $K\in \cK_o^n$ such that  \begin{eqnarray*} 
	\frac{\int_{\sphere} g(v)\,d\mu(v) }{|\mu|} = \frac{ \int_{S^{n-1}}g(v)\,d\cV(K,v)}{ \int_{S^{n-1}}\,d\cV(K,v)}=\frac{\int_{S^{n-1}}g(\alpha_K(u))\hat\varphi(\rho_K(u))\phi_{2}(u)\,du}{\int_{S^{n-1}} \hat\varphi(\rho_K(u))\phi_{2}(u)\,du}
	 \end{eqnarray*}
	hold for each bounded Borel function $g: S^{n-1}\rightarrow \R$. \end{itemize}  \ec
	
\begin{proof} As explained in Section \ref{Section:2}, under the conditions given in Corollary \ref{cor-5-1}, $\phi(x)=\psi(|x|)\phi_{2}(\bar{x})$ satisfies conditions C1) and C2). The argument in ii) is equivalent to $$\frac{\mu}{|\mu|}=\frac{\widetilde{C}_{\phi,\mathcal{V}}(K,\cdot)}{\widetilde{C}_{\phi,\mathcal{V}}(K,S^{n-1})}. $$ The equivalence between i) and ii) is an immediate consequence from   (\ref{concentration-0--10}), (\ref{concentration-1--1}), Corollary \ref{measure-change}  and Theorem \ref{solution-general-dual-Orlicz-main theorem}.\end{proof}

 \section{Uniqueness of solutions of the general dual Orlicz-Minkowski problem} \label{Section:6}
	
It seems very difficult and maybe even impossible to obtain the uniqueness of solutions of the general dual Orlicz-Minkowski problem for general $\phi$. In this section, the uniqueness will be proved in special cases. In order to get this done, we need the following theorem. 
 
 \bt\label{uniqueness-1-1} Let $\phi$ be a function satisfying condition C1)  and that $\phi(x)|x|^n$ is strictly radially decreasing on $\Rn\setminus\{o\}$. If $K,L\in \cK_o^n$ satisfy that $\cV(K,\cdot)=\cV(L,\cdot)$, then $K=L$. \et

The proof of Theorem \ref{uniqueness-1-1} follows an argument similar to those in \cite{zhao, ZSY2017}, and heavily relies on  \cite[Lemma 5.1]{zhao}. For readers' convenience, we list  \cite[Lemma 5.1]{zhao} below as Lemma \ref{uniqueness for dual quermassintegral} and provide a brief sketch of the proof of Theorem \ref{uniqueness-1-1}.  
	\bl\label{uniqueness for dual quermassintegral}
	Suppose that $K', L\in\mathcal{K}_{o}^n$. If the following sets 
	\begin{eqnarray*}
		\eta_1&=&\big\{v\in S^{n-1}: \ \ h_{K'}(v)>h_{L}(v) \big\}, \\
		\eta_2 &=& \big\{v\in S^{n-1}: \ \ h_{K'}(v)<h_{L}(v) \big\}, \\ 
		\eta_{3}&=&\big \{v\in S^{n-1}: \ \ h_{K'}(v)=h_{L}(v)\big\}
	\end{eqnarray*} are nonempty, then the following statements are true:
	\begin{itemize} 
	\item[a)] if $u\in\pmb{\alpha}_{K'}^*(\eta_1)$, then $\rho_{K'}(u)>\rho_{L}(u)$;
	\item[b)] if $u\in\pmb{\alpha}_{L}^*(\eta_2\cup\eta_3)$, then $\rho_{L}(u)\geq \rho_{K'}(u)$;
	\item[c)]  $\pmb{\alpha}_{K'}^*(\eta_1)\subset\pmb{\alpha}_{L}^*(\eta_1)$;
	\item[d)] $\mathcal{H}^{n-1}(\pmb{\alpha}_{L}^*(\eta_1))>0\ \ \mbox{and}\ \ \mathcal{H}^{n-1}(\pmb{\alpha}_{K'}^*(\eta_2))>0$.
	\end{itemize}
	\el
	
\vskip 2mm \noindent {\em Proof of Theorem \ref{uniqueness-1-1}.} 
		Assume that $K, L\in \cK_o^n$ with $\cV(K,\cdot)=\cV(L,\cdot)$  are not dilates of each other, namely, $K\neq t L$ for any $t>0$. Hence, there exists some constant $t_0>0$ such that  $K'=t_{0}K$ is a convex body with $\eta_1,\eta_2,\eta_3$ defined in Lemma \ref{uniqueness for dual quermassintegral} being nonempty.   
		
		Recall that $\Psi_{K}(u)=\phi(\rho_K(u)u)[\rho_K(u)]^n$ for $u\in\sphere.$ Due to Lemma \ref{uniqueness for dual quermassintegral}  and the fact that  $\phi(x)|x|^n$ is strictly radially decreasing  on $\Rn\setminus\{o\},$ one has, for all $u\in \pmb{\alpha}^{*}_{K'}(\eta_{1})$, \begin{equation}0< \Psi_{K'}(u)=\phi(\rho_{K'}(u)u)[\rho_{K'}(u)]^n < \phi(\rho_{L}(u)u)[\rho_{L}(u)]^n=\Psi_{L}(u). \label{comparison--1-1} \end{equation}    Now we claim that  the spherical measure of $\pmb{\alpha}^{*}_{K'}(\eta_{1})$ is positive. In fact, this claim follows from  Definition \ref{general dual Orlicz curvature measure} and  Lemma \ref{uniqueness for dual quermassintegral}  as follows:  $$\int_{\pmb{\alpha}^{*}_{K}(\eta_{1})}\Psi_{K}(u)\,du=\cV(K,\eta_{1}) = \cV(L,\eta_{1})=\int_{\pmb{\alpha}^{*}_{L}(\eta_{1})}\Psi_{L}(u)\,du  >0.$$ Moreover, by  (\ref{comparison--1-1}) and  Lemma \ref{uniqueness for dual quermassintegral}, one has   \begin{equation*} \cV(K,\eta_1) =\int_{\pmb{\alpha}^{*}_{L}(\eta_{1})}\Psi_{L}(u)\,du \geq  \int_{\pmb{\alpha}^{*}_{K'}(\eta_{1})}\Psi_{L}(u)\,du> \int_{\pmb{\alpha}^{*}_{K'}(\eta_{1})}\Psi_{K'}(u)\,du>0.\end{equation*}  Due to the easily checked fact $\pmb{\alpha}^{*}_{K'}(\eta_{1})=\pmb{\alpha}^{*}_{K}(\eta_{1})$ and Definition \ref{general dual Orlicz curvature measure}, one gets \begin{eqnarray*} \cV(K,\eta_1)&=& \int_{\pmb{\alpha}^{*}_{K}(\eta_{1})}\Psi_{K}(u)\,du \\ &=& \int_{\pmb{\alpha}^{*}_{K'}(\eta_{1})}\phi(\rho_K(u)u)[\rho_K(u)]^n\,du \\ &>&  \int_{\pmb{\alpha}^{*}_{K'}(\eta_{1})}\Psi_{K'}(u)\,du\\ &=& \int_{\pmb{\alpha}^{*}_{K'}(\eta_{1})} \phi(t_0\rho_K(u)u)[t_0\rho_K(u)]^n \,du >0.\end{eqnarray*}  Together with the fact that $\phi(x)|x|^n$ is strictly radially decreasing on $\Rn\setminus\{o\}$, one has  $t_0>1$ and moreover  \begin{equation} \label{comparison-0--1-0} \phi(\rho_K(u)u)[\rho_K(u)]^n>\phi(t_0\rho_K(u)u)[t_0\rho_K(u)]^n \end{equation}  holds for all $u\in \sphere$. 
		
	Similarly, one can check that the spherical measure of $\pmb{\alpha}_{L}^{*}(\eta_2)$ is positive. It follows from Lemma \ref{uniqueness for dual quermassintegral} that  $\pmb{\alpha}_{L}^{*}(\eta_2)\subseteq \pmb{\alpha}_{K'}^{*}(\eta_2)$ and  
		\[0< 
		\cV(K,\eta_2)=\cV(L,\eta_2)= \int_{\pmb{\alpha}_{L}^{*}(\eta_2)}\Psi_{L}(u)\,du\leq \int_{\pmb{\alpha}_{K'}^{*}(\eta_2)}\Psi_{K'}(u)\,du=\cV(K',\eta_2).
		\] Together with (\ref{comparison-0--1-0}), Definition \ref{general dual Orlicz curvature measure}, and $\pmb{\alpha}^{*}_{K'}(\eta_{2})=\pmb{\alpha}^{*}_{K}(\eta_{2})$, one has $$\cV(K,\eta_2)\leq\cV(K',\eta_2)<\cV(K,\eta_2).$$  This is impossible, and hence $K$ and $L$ are dilates of each other. 
		
		Now we claim that $K=L$. Assume not, i.e., there exists a constant $t\neq 1$ such that $K=tL$. Let $t>1$ and hence $\phi(\rho_L(u)u)[\rho_L(u)]^n >\phi(\rho_{ K}(u)u)[\rho_{ K}(u)]^n$ for all $u\in \sphere$. We can get a contradiction as follows:  \begin{eqnarray*}
			\cV(K,\sphere)&=&\cV(L,\sphere)\\
			&=& \int_{\sphere} \phi(\rho_L(u)u)[\rho_L(u)]^n\,du\\ 
			&>& \int_{\sphere}\phi(\rho_{ K}(u)u)[\rho_{ K}(u)]^n\,du\\
			&=&\cV(K,\sphere),
		\end{eqnarray*}  where we have used the assumption that $\cV(K,\cdot)=\cV(L,\cdot)$.
		
		 Similarly, one can show that $t<1$ is not possible, and hence  $K=L$ as desired. \hfill $\Box$
	
\vskip 2mm \noindent {\bf Remark.} When $\phi(x)=\psi(|x|)\phi_{2}(\bar{x})$ as stated in Case 2 in Section \ref{Section:2}, then $\phi(x)|x|^n$ is a strictly radially decreasing function if  $\hat{\varphi}(t)=  n\psi(t)t^n$ is a strictly decreasing function on $t\in (0, \infty)$. For instance, if $\phi(x)=\|x\|^{q-n}$ for $q<0$, then $$\phi(x)|x|^n=\|x\|^{q-n}|x|^n=|x|^q\|\bar{x}\|^{q-n}$$  is a  strictly radially decreasing function on $\Rn\setminus\{o\}$.  On the other hand, if $\phi$ is smooth enough, say the gradient of $\phi$ (denoted by $\nabla \phi$) exists on $\Rn\setminus \{o\}$, 	a typical condition to make $\phi(x)|x|^n$ strictly radially decreasing is  $\big\langle \nabla\big(\phi(x)|x|^n\big), x\big\rangle <0$  or equivalently   $\langle \nabla \phi(x), x\rangle +n\phi(x)<0$   for all $x\in\Rn\setminus \{o\}.$  
	
	We are now ready to state our result regarding the uniqueness of solutions to the general dual Orlicz-Minkowski problem. If $\phi_{2}(u)=1$ for all $u\in \sphere$, it goes back to the case proved by Zhao \cite{zhao}.  \bc \label{Cor-6-1}
	Let $\phi(x)=|x|^{q-n} \phi_{2}(\bar{x})$ with $q<0$ and $\phi_{2}: \sphere\rightarrow (0, \infty)$ a positive continuous function. Then the following statements are equivalent: \begin{itemize} \item [i)] $\mu$ is a nonzero finite Borel measure on $S^{n-1}$ satisfying (\ref{condition for Minkowski problem}); 
		\item [ii)] there exists a unique convex body $K\in\cK_o^n$, such that, $\mu=\widetilde{C}_{\phi,\mathcal{V}}(K,\cdot)$. \end{itemize} \ec
		\begin{proof} The argument from ii) to i) follows along the same lines as the arguments for  (\ref{concentration-0--10}) and (\ref{concentration-1--1}). 
		On the other hand, it follows from Theorem \ref{solution-general-dual-Orlicz-main theorem} that, if  $\mu$ is a nonzero finite Borel measure on $S^{n-1}$ satisfying (\ref{condition for Minkowski problem}), then  there is a convex body $\widetilde{K}\in\mathcal{K}_{o}^n$ such that
		$$\frac{\mu}{|\mu|}=\frac{\widetilde{C}_{\phi,\mathcal{V}}(\widetilde{K},\cdot)}{\widetilde{C}_{\phi,\mathcal{V}}(\widetilde{K},S^{n-1})}. $$ 
				By Corollary \ref{measure-change}, $\pmb{\alpha}^*_{\lambda K}(\eta)=\pmb{\alpha}^*_{K}(\eta)$ and $\rho_{\lambda  K}=\lambda \rho_K$  for any constant $\lambda>0$, and the fact that  $u\in \pmb{\alpha}_K^*(\eta)$  if and only if $\alpha_K(u)\in \eta$ (see  \cite[(2.21)]{HLYZ}), one has,  for any $\lambda>0$ and for any Borel set $\eta\subseteq \sphere$, 
		\begin{eqnarray}
			\cV(\lambda K,\eta) &=&\int_{\pmb{\alpha}^*_ {\lambda K}(\eta)} \big[\rho_{\lambda  K}(u)\big]^{q}\phi_{2}(u)du \nonumber\\
			&=&\lambda^{q} \int_{\pmb{\alpha}^*_ K(\eta)}\big[\rho_{K}(u)\big]^{q}\phi_{2}(u)du \nonumber\\
			&=&\lambda^{q}\cV(K,\eta). \label{homogenes-3-3}
		\end{eqnarray} Hence, $\cV(\lambda K,\cdot) =\lambda^{q}\cV(K,\cdot)$ and  
		\[
		\mu=\frac{|\mu|}{\widetilde{C}_{\phi,\mathcal{V}}(\widetilde{K},\sphere)}\widetilde{C}_{\phi,\mathcal{V}}(\widetilde{K},\cdot)=\widetilde{C}_{\phi,\mathcal{V}}(K,\cdot),
		\] where
		\[
		K=\bigg(\frac{|\mu|}{\widetilde{C}_{\phi,\mathcal{V}}(\widetilde{K},\sphere)}\bigg)^{\frac{1}{q}}\widetilde{K}.
		\] Hence, $K\in \cK_o^n$ is a convex body such that $\mu=\cV(K, \cdot)$, if $\mu$ is a nonzero finite Borel measure on $S^{n-1}$ satisfying (\ref{condition for Minkowski problem}).  The uniqueness of $K$ is an immediate consequence of Theorem \ref{uniqueness-1-1} and the remark after its proof.  
	\end{proof}
	
	The solution for $\mu$ being a discrete measure is stated in the following proposition.

	\bp Let $\phi(x)=|x|^{q-n} \phi_{2}(\bar{x})$ with $q<0$ and $\phi_{2}: \sphere\rightarrow (0, \infty)$ a positive continuous function. Suppose that $\mu=\sum_{i=1}^m \lambda_i\delta_{u_i}$ with all $\lambda_i>0$ is a discrete measure not concentrated in any closed hemisphere  (i.e., satisfying (\ref{condition for Minkowski problem})). Then, there exists a unique polytope $P\in\cK_o^n$, such that, $\mu=\widetilde{C}_{\phi,\mathcal{V}}(P,\cdot)$ and $u_1,u_2,\cdots, u_m$ are the unit normal vectors of the faces of $P$. \ep
	\begin{proof}  It follows from Corollary \ref{Cor-6-1} that there exists a unique convex body $K_0\in \cK_o^n$, such that, $\mu=\widetilde{C}_{\phi,\mathcal{V}}(K_0,\cdot)$. The desired argument in this proposition follows if we can prove that $K_0$ is a polytope with $u_1,u_2,\cdots, u_m$ being the unit normal vectors of its faces. To this end, let $M\in \cK_o^n$ be a polytope circumscribed about $K_0$ whose faces have  the unit normal vectors being exactly $u_1,u_2,\cdots, u_m.$ Hence $K_0\subseteq M$ and $h_M(u_i)=h_{K_0}(u_i)$ for all $i=1, 2, \cdots, m$.

	 Suppose that $K_0\neq M$ (as otherwise, nothing to prove). In this case, there exists a set $\eta_M\subseteq \sphere,$ such that, the spherical measure of $\eta_M$ is positive and $\rho_M(u)>\rho_{K_0}(u)$ on $\eta_M$. It follows from (\ref{homogenes-1-1}) and (\ref{homogenes-2-2}) that $\eV(M)<\eV(K_0)$  and $\cV(L, \sphere)=-q\eV(L)$ for all $L\in \cK_o^n$. Hence, $\cV(M, \sphere)<\cV(K_0,\sphere)=|\mu|$. By (\ref{homogenes-3-3}), there exists a constant $0<c<1$, such that 
	$$\cV(cM, \sphere)=\cV(K_0,\sphere)=|\mu|.$$ On the other hand, from Corollary \ref{Cor-6-1} and  the proof of Theorem \ref{solution-general-dual-Orlicz-main theorem},  the convex body $(-q)^{1/q}K_0\in \cK_o^n$ is the unique convex body  such that $\eV\big((-q)^{1/q}K_0\big)=|\mu|$ and 
	 $$\mathcal{F}\big((-q)^{1/q}K_0\big)=\sup\big\{\mathcal{F}(K): \eV(K)=|\mu| \ \mathrm{and}\ K\in \cK_o^n \big\}.$$ However, this is impossible because $\eV\big((-q)^{1/q}c M\big)=|\mu|$ and 
	 \begin{eqnarray*} \mathcal{F}\big((-q)^{1/q}c M\big)&=&-\frac{1}{|\mu|}\int_{S^{n-1}}\big[\log h_M(v)+\log c+ \log (-q)/q\big]\,d\mu(v) \\&>&-\frac{1}{|\mu|}\int_{S^{n-1}}\big[\log h_M(v)+ \log (-q)/q\big]\,d\mu(v)\\ &=& -\frac{1}{|\mu|} \cdot \sum_{i=1}^m  \lambda_i \big[\log h_M(u_i)+ \log (-q)/q\big] \\ &=& -\frac{1}{|\mu|} \cdot \sum_{i=1}^m  \lambda_i \big[\log h_{K_0}(u_i)+ \log (-q)/q\big] \\ &=& \mathcal{F}((-q)^{1/q}K_0), \end{eqnarray*} where the inequality is due to $0<c<1$. Hence $M=K_0$ is a polytope.  Moreover, it is easy to get the relation between $\lambda_i$ and the polytope $K_0$. In fact, \begin{eqnarray*} \lambda_i&=&\int_{\{u_i\}}\,d\mu =\int_{\{u_i\}}\,d\cV(K_0, v) \\ &=& \int_{\pmb{\alpha}^*_ {K_0}(u_i)} \big[\rho_{K_0}(v)\big]^{q}\phi_{2}(v)\,dv\\&=& \int_{\nu_ {K_0}^{-1}(\{u_i\})}\langle x,\nu_{K_0}(x)\rangle  |x|^{q-n}\phi_{2}(\bar{x}) \,d\mathcal{H}^{n-1}(x)\\  &=& \int_{\nu_ {K_0}^{-1}(\{u_i\})}\langle x,\nu_{K_0}(x)\rangle  \phi(x) \,d\mathcal{H}^{n-1}(x)>0\end{eqnarray*} where the third equality follows from (\ref{homogenes-3-3}) and the fourth equality follows from Corollary \ref{measure-change}. Let $P=K_0$, and then $P$ is the desired polytope, such that, $\mu=\widetilde{C}_{\phi,\mathcal{V}}(P,\cdot)$ and $u_1,u_2,\cdots, u_m$ are the unit normal vectors of the faces of $P$.
	\end{proof}

	\vskip 2mm \noindent {\bf Acknowledgments.}  The research of DY  is supported by a NSERC
grant. The authors are greatly indebted to the referee for many valuable comments which improve largely the quality of the paper.

	\vskip 2mm \noindent Sudan Xing, \ \ \ {\small \tt sudanxing@gmail.com}\\
{ \em Department of Mathematics and Statistics,   Memorial University of Newfoundland,
   St.\ John's, Newfoundland, Canada A1C 5S7 }

\vskip 2mm \noindent Deping Ye, \ \ \ {\small \tt deping.ye@mun.ca}\\
{ \em Department of Mathematics and Statistics,
   Memorial University of Newfoundland,
   St.\ John's, Newfoundland, Canada A1C 5S7 }

\end{document}